%% file: Nature_R0_v5_submitted.tex
\documentclass[pdflatex,sn-mathphys-num]{sn-jnl}

\usepackage{graphicx}%
\usepackage{multirow}%
\usepackage{amsmath,amssymb,amsfonts}%
\usepackage{amsthm}%
\usepackage{mathrsfs}%
\usepackage[title]{appendix}%
\usepackage{xcolor}%
\usepackage{textcomp}%
\usepackage{manyfoot}%
\usepackage{booktabs}%
\usepackage{pdfpages}
\usepackage{algorithm}%
\usepackage{algorithmicx}%
\usepackage{algpseudocode}%
\usepackage{listings}%
\usepackage{comment}%

\theoremstyle{thmstyleone}
\theoremstyle{thmstyletwo}

\theoremstyle{thmstylethree}

\usepackage{tabularx}
\usepackage{makecell}
\usepackage{rotating}

\definecolor{darkgreen}{rgb}{0.0, 0.5, 0.0}

\usepackage{verbatim}

\usepackage{rotating}
\usepackage{float}
\usepackage{placeins}

\usepackage{graphicx}
\usepackage{subcaption}

\begin{document}

\title[Article Title]{State-wise Economic Viability of Long-Duration Energy Storage Systems in the United States}


\author*[1]{\fnm{Alexandre} \sur{Moreira}}\email{AMoreira@lbl.gov}
\equalcont{These authors contributed equally to this work.}

\author[2]{\fnm{Patricia} \sur{Silva}}\email{patriciadesousaoliveira@gmail.com}
\equalcont{These authors contributed equally to this work.}

\author[1]{\fnm{Miguel} \sur{Heleno}}\email{MiguelHeleno@lbl.gov}

\author[2]{\fnm{André} \sur{Marcato}}\email{amarcato@ieee.org}

\affil[1]{\orgdiv{Energy Storage and Distributed Resources Division}, \orgname{Lawrence Berkeley National Laboratory}, \orgaddress{\street{1 Cyclotron Rd}, \city{Berkeley}, \postcode{94720-8134}, \state{California}, \country{United States}}}

\affil[2]{\orgdiv{Departament of Electrical Engineering}, \orgname{Federal University of Juiz de Fora}, \orgaddress{\street{Rua José Lourenço Kelmer}, \city{Juiz de Fora}, \postcode{36036-330}, \state{Minas Gerais}, \country{Brazil}}}


\abstract{Long-duration energy storage (LDES) assets can be fundamental resources for the next-generation power systems. However, LDES technologies are still immature and their future technology costs remain highly uncertain. In this context, we perform in this paper an extensive study to estimate the maximum LDES technology costs (which we define as viability costs) under which LDES systems would be economically viable in each state of the contiguous U.S. according to their characteristics.  Our results indicate that only 4 states (out of 48) would be able to remove firm conventional generation supported by LDES systems without increasing their total system costs under the current US-DOE cost target of US\$ 1,100 $\text{kW}^{-1}$ for multi-day LDES. In addition, we find that states with the highest LDES viability costs have in general low participation of thermal generation, a high share of wind generation, and higher thermal-related fixed operation and maintenance (FO\&M) costs.
}




\maketitle

\section{Introduction}\label{sec:introduction}

Long-duration energy storage (LDES) systems can play a pivotal role in future power systems as they have the capability to charge and discharge over extended periods of time either leveraging excess or compensating for insufficient momentaneous generation capacity. This capability allows LDES technologies to counterbalance a potential lack of sufficient availability of intermittent generation capacity for a sustained amount of time \cite{Zhang2020}, which is a service that short-duration energy storage (SDES) systems cannot support. 
Nonetheless, in the heterogeneous landscape of the energy systems in the United States, it is sill unclear what the thresholds would be for the economic viability of LDES systems in the regional power systems of the future.

In general, the economic viability of LDES technologies is assessed in the current literature by means of assuming future technology cost projections. For instance, \cite{colombo2023value} focuses on a capacity expansion study for California targeting years 2030, 2035, 2040, and 2045 while taking into account, as candidate resources, 100-h LDES systems with 50\% and 80\% round-trip efficiency. The LDES investment costs are assumed to range from US\$ 1,715 $\text{kW}^{-1}$ in 2030 to US\$ 1,500 $\text{kW}^{-1}$ and, based on this assumption, the authors in \cite{colombo2023value} conclude that only the 80\% round trip efficiency (RTE) LDES would be competitive with 4-h lithium-ion batteries. In \cite{staadecker2024value}, the authors examine how factors such as the share of wind and solar generation capacity, the availability of hydropower generation, and the deployment of new transmission capacity can impact the appetite for investments in LDES. Their study combines these factors into 39 capacity expansion scenarios, using the Western Electricity Coordinating Council (WECC) system and  assuming, for candidate storage resources, overnight power capacity costs of US\$ 19.58 $\text{kW}^{-1}$, along with annual operation and maintenance costs of US\$ 6.10 $\text{kW}^{-1}$, and baseline energy capacity cost of US\$ 22.43 $\text{kWh}^{-1}$. The Study finds that (i) wind-dominant regions would be ideally served by 10-h to 20-h storage resources, (ii) a reduced hydropower availability in hydro-dominant areas would increase the ideal duration of newly installed storage assets from around 6-h to 23-h, (iii) a variation of overnight energy capacity costs would imply in ideal durations ranging from 9-h to around 800-h. In \cite{guerra2020}, the economic viability of different types of storage is assessed for the Western US power system by a proposed 3-step approach that consists of obtaining a capacity expansion solution, running a production cost analysis without seasonal storage and then repeating this run with seasonal storage. Based on their results, considering a given cost range, hydrogen seasonal storage with 1-week duration would be economically viable with power and energy technology costs lower or equal to around US\$ 1507 $\text{kW}^{-1}$ and US\$ 1.8 $\text{kWh}^{-1}$, respectively, by 2025. In \cite{hunter2021techno}, the capacity expansion and production cost model results from \cite{Zhang2020} are leveraged to compare the cost-effectiveness of different LDES technologies, in terms of levelized cost of energy (LCOE), to flexible power generation systems, such as natural gas combined cycle plants (NG-CC). Based on the considered reference technology costs, NG-CC are more competitive than any other options while, among the storage technologies, diabatic compressed air energy storage in a salt cavern (D-CAES|Salt) achieves the lowest LCOE values ranging from lower than US\$ 200 MWh$^{-1}$ for durations between 12 and 72 h to lower than US\$ 400 MWh$^{-1}$ for durations between 72 to 168 h (in 2018 USD). However, these values can dramatically change due to uncertainty in future investment costs of these technologies. Furthermore, \cite{sepulveda2021design} indicates that (i) power systems dominated by intermittent resources can have more than 10\% decrease in total costs if the LDES energy capacity costs are limited to US\$ 20 $\text{kWh}^{-1}$, (ii) LDES energy capacity costs would have to fall lower than US\$ 10 $\text{kWh}^{-1}$ to displace nuclear plants and lower than US\$ 1 $\text{kWh}^{-1}$ to able to displace natural gas with carbon capture storage (CCS).         

The relevance of the aforementioned works notwithstanding, emerging LDES technologies are currently immature and their future cost projections are highly uncertain. In this context, we introduced in \cite{first_paper} two capacity expansion-based models to compute the viability costs below which LDES technologies would be economically viable from the system's perspective. The first one is a {\it baseline model} where firm flexible generators are considered without LDES and the second one is an {\it opportunity value model} in which conventional flexible generators are removed and a given capacity (which varies through a sensitivity analysis) of LDES is included. To compute the LDES viability costs, we first obtain the maximum opportunity value for LDES given its main parameters (duration and efficiency) and the energy matrix of the system under consideration. We define this LDES maximum opportunity value as the maximum avoided cost that can be achieved when replacing conventional firm generators with a mix of intermittent power plants along with short-duration energy storage (SDES) and LDES systems. Then, the LDES viability costs are calculated as ratio between the LDES maximum opportunity value and its corresponding capacity.   

In this paper, we expand the methodology proposed in \cite{first_paper} to calculate the LDES viability costs across the 48 contiguous U.S. states, focusing on 100-h duration systems, and to conduct an extensive analysis of the economic viability of LDES. Our analysis is based on information from comprehensive datasets from the NREL's Cambium 2022 \cite{gagnon2023cambium_site}, the Regional Energy Deployment System (ReEDS) base \cite{reedsGitHub}, the 2022 Annual Technology Baseline (ATB)  data \cite{atb2022}, and the Annual Energy Outlook (AEO) 2023 website information \cite{AEO_2023}. Our results reveal a wide range of viability costs for LDES technologies across the States, from -US\$ 11.81 $\text{kW}^{-1}$ to US\$ 5,993.94 $\text{kW}^{-1}$. While the US-DOE has set a 2030 cost target of US\$ 1,100 $\text{kW}^{-1}$ for multi-day LDES (36–160 hours) \cite{Liftoff_DOE_2023}, our results indicate that this target would only be economically viable in 4 states (KS, ND, NE, and VA). Lowering the target to US\$ 500 $\text{kW}^{-1}$ by 2050 would increase this number of states to 9 (including CA, CO, ID, MT, and NM). A further reduction to US\$ 300 $\text{kW}^{-1}$ would make LDES viable in 17 states, expanding the list to include TX, WA, and MN, among others. Our analysis also identified five states with negative viability costs for all considered LDES power capacities: Alabama (AL), Connecticut (CT), Delaware (DE), New Jersey (NJ), and  Ohio (OH). These negative values indicate that total system costs would increase if conventional firm generators are replaced with intermittent plants and storage systems (SDES and LDES) and therefore LDES is not viable. Moreover, while the Pathways to Commercial Liftoff report \cite{Liftoff_DOE_2023} estimates a need for 225–460 GW of LDES by 2050, our findings suggest that, to replace firm generators, the minimum required capacity of LDES systems (100-hour duration) is approximately 646.09 GW (64.61 TWh) for the 43 states where maximum viability costs are positive, with 39.53\% of these 43 states requiring more than 10 GW (1 TWh).

Furthermore, our results also identify system characteristics that can strongly influence the economic viability of LDES. For instance, states with the highest viability cost tend to have low participation of thermal (referring to gas and coal throughout this paper) generation (with thermal capacity utilization often below 20\%), a high share of wind generation in the intermittent mix (above 65\%), and higher fixed operation and maintenance (FO\&M) costs related to thermal plants. In contrast, states with the lowest LDES viability costs usually have greater reliance on thermal generation, lower average hourly capacity factor (CF) of intermittent energy sources (IES)  compared to the national average (which is around 27\% considering 48 states) and a high share of solar generation in the intermittent portfolio ($>65\%$).

\section{Results}\label{sec:results}

\subsection{Wide range of viability costs and low ambition in LDES technology targets}\label{sec:map_boundary_cost}

\begin{figure}[!t]
    \centering
     \includegraphics[width=1.0\textwidth,height=1.0
     \textheight,keepaspectratio]{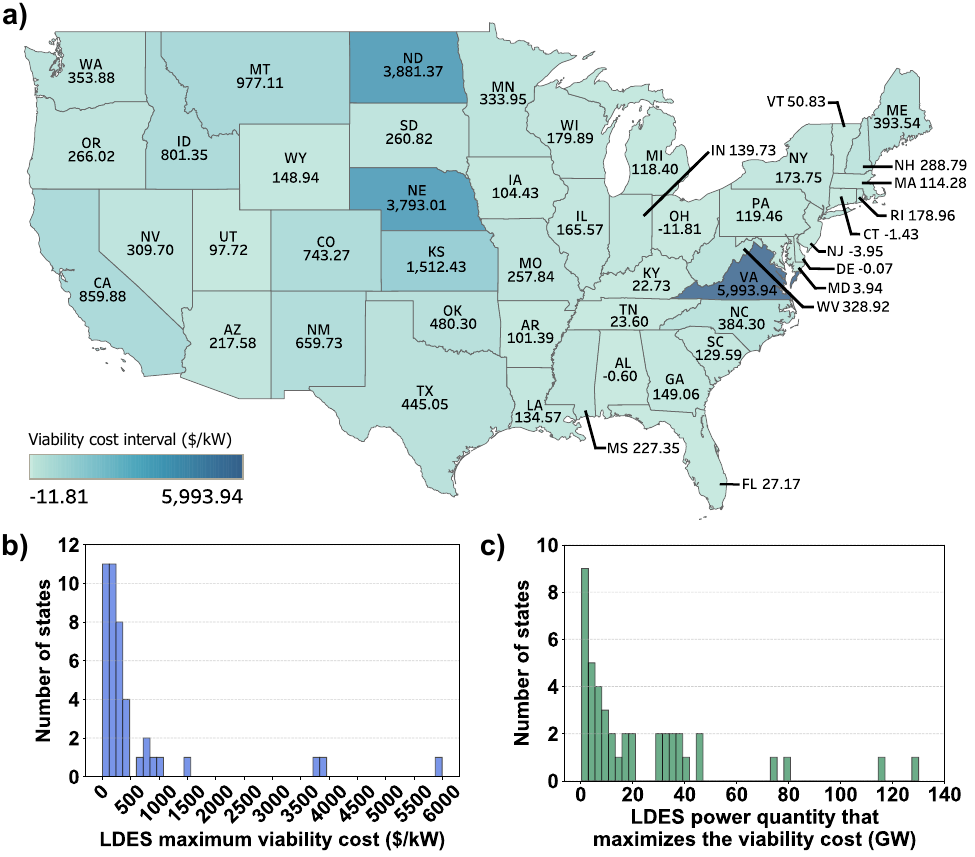}
     \caption{LDES maximum viability cost and power capacity distribution across U.S. states in 2050. (a) Maximum viability cost of 100-hour LDES (\$/kW) for each U.S. state, considering the power capacity that maximizes the viability cost. (b) Histogram of the LDES maximum viability cost (\$/kW) across states with positive values. (c) Histogram of the LDES power capacity (GW) that maximizes the viability cost in states with positive values.}
    \label{Fig:Map}
\end{figure}

Leveraging the optimization models proposed in \cite{first_paper}, we estimate here the maximum LDES viability cost for each state of the contiguous U.S., considering their 2050 energy matrices according to NREL's Cambium 2022 dataset \cite{gagnon2023cambium_site}, and depict the resulting values in Fig. \ref{Fig:Map}. Our study evaluates these maximum LDES viability costs based on replacing gas and coal generators with a combination of intermittent units alongside SDES and LDES resources and computing the maximum avoided costs (with respect to the total costs of the original system) - hereinafter referred to as maximum opportunity values - that can be achieved by various power capacities of a 100-hour LDES system with 42.5\% RTE (realistic parameters of an iron-air storage \cite{DOE_LDES_iron_air_demonstration_2024,McKerracher2015_IronAirReview}). More specifically, for each state, we consider a range of 100-h LDES power capacities varying from 0 to 150 GW and calculate the corresponding LDES viability cost as the ratio between the resulting maximum opportunity value and the considered capacity. As a result, each state will have a certain 100-h LDES power capacity that maximizes its LDES viability cost.

\begin{figure}[!b]
    \centering
     \includegraphics[width=1.0\textwidth,height=1.0
     \textheight,keepaspectratio]{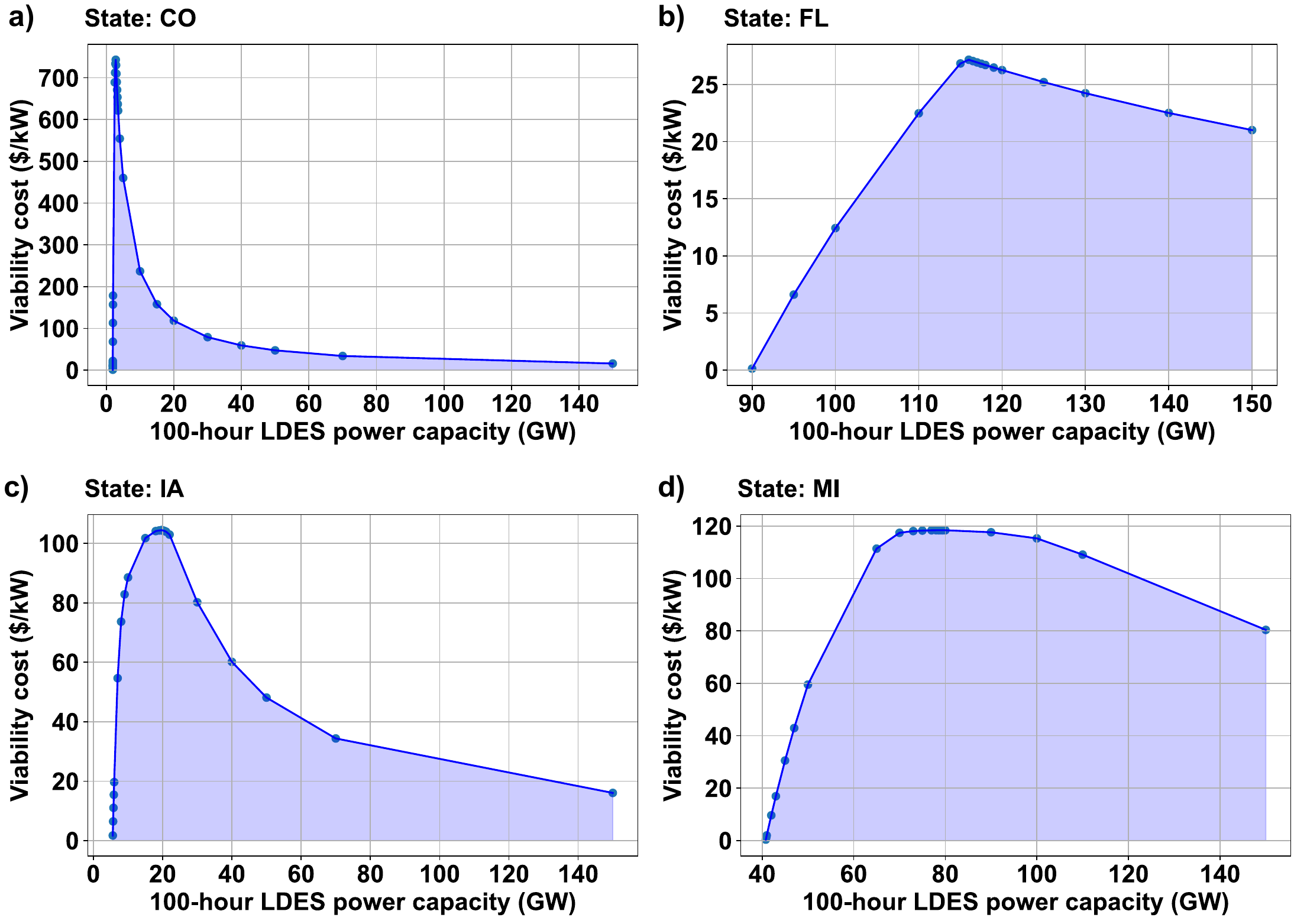}
    \caption{The viability costs of LDES below which these technologies will be economically viable for the system in 2050, for: (a) Colorado, (b) Florida, (c) Iowa, and (d) Michigan.}    \label{Fig:boundary_cost_of_4_states}
\end{figure}

The maximum LDES viability cost across U.S. states varies significantly, ranging from -US\$ 11.81 $\text{kW}^{-1}$ to US\$ 5,993.94 $\text{kW}^{-1}$. In states where the viability cost is negative, namely AL, NJ, OH, CT, and DE, 100-h LDES is not viable, since the total system costs increase when thermal generation is replaced by intermittent plants. These states tend to have limited availability of intermittent resources, which are mostly solar. More specifically, 4 of these states have an average hourly CF for intermittent units below 25\% and also 4 of them have solar capacity exceeding 95\% of their intermittent portfolio, which limits the flexibility of charging and discharging cycles. In addition, all of them rely on thermal generation to supply more than 30.00\% of their annual demand and 4 of them incur in thermal-related FO\&M costs that represent less than 10.00\% to overall costs, indicating dependence on low cost firm generators. The dependence on these assets reduces the benefits of retiring thermal assets, ultimately limiting the opportunity value of LDES in these regions.

Surprisingly, according to our results, 100-h LDES systems would be economically viable for thermal replacement in only 4 states (KS, ND, NE, and VA) under the US-DOE's 2030 technology cost target of US\$ 1,100 $\text{kW}^{-1}$ for multi-day LDES \cite{Liftoff_DOE_2023}. If the US-DOE target falls to US\$ 500 $\text{kW}^{-1}$ by 2050, 9 states would surpass this threshold, including CA, CO, ID, KS, ND, NE, NM, MT, and VA. A further reduction to US\$ 300 $\text{kW}^{-1}$ would result in 17 states exceeding this value, expanding the list to include ME, MN, NC, NV, OK, TX, WA, and WV.

\subsection{LDES capacity required for thermal replacement}\label{sec:power_quantity_required}

\begin{figure}[!tb]
    \centering
     \includegraphics[width=1.0\textwidth,height=1.0
     \textheight,keepaspectratio]{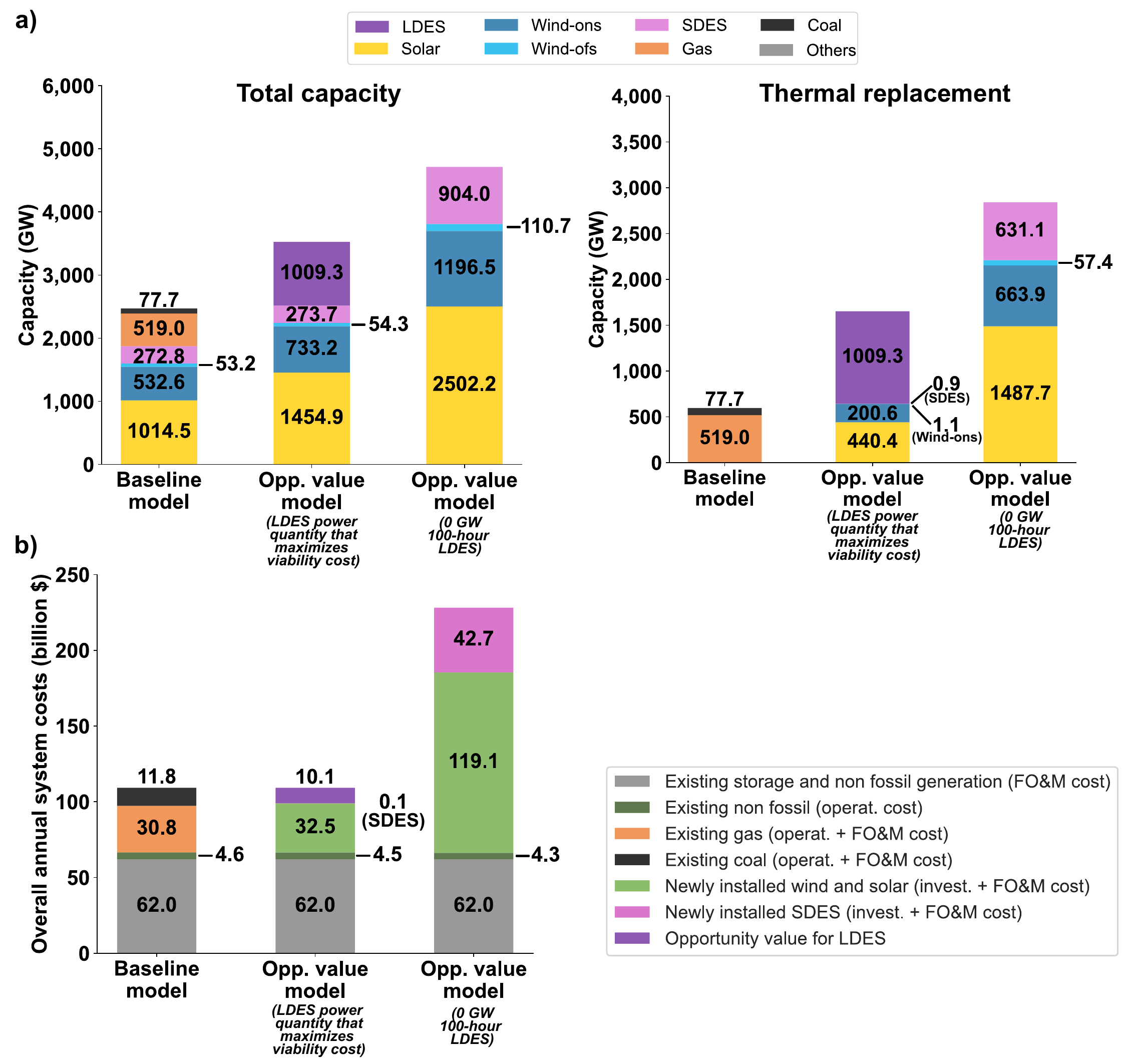}
    \caption{(a) U.S. total system capacity by technology in 2050: \textit{baseline model} vs. \textit{opportunity value model} with and without LDES. The first bar represents the total system capacity when considering gas and coal power plants with no inclusion of LDES. The second bar shows the scenario where gas and coal power plants are replaced by a combination of LDES, SDES, and intermittent energy sources, with LDES capacity set to the amount that maximizes the viability cost. The third bar represents a scenario where gas and coal power plants are also replaced by SDES and intermittent energy sources, but without the inclusion of LDES. The three bars on the right show the thermal replacement. The first bar indicates the amount of gas and coal capacity in the baseline, while the next two show the net increase in other technologies required to replace the output previously provided by gas and coal. (b) Same analysis as in (a), but for cost values instead of capacity. Both figures include data from 43 states, excluding values those where the maximum viability cost is negative.}
    \label{Fig:National_capacity_baseline_and_opp_value}
\end{figure}

At a national scale, our analysis reveals a substantial gap between projected LDES needs and the minimum LDES capacities that would be required to retire thermal generators in the United States while maintaining the same overall system costs. While the Pathways to Commercial Liftoff report \cite{Liftoff_DOE_2023} suggests that the U.S. grid may need between 225 and 460 GW of LDES capacity in 2050, including multi-day/week LDES (36–160 hours) as well as inter-day storage (10–36 hours), our results estimate that the minimum required 100-hour LDES capacity is approximately 646.09 GW (64.61 TWh) for the 43 states where maximum viability costs are positive. More specifically, 39.53\% of these 43 states would need more than 10 GW, while only 16.28\% would require less than 1 GW. Considering the 100-h LDES power capacities that maximize the viability costs, these numbers increase significantly to 1,009.30 GW (100.93 TWh) for these 43 states. In this case, 53.49\% of these states are associated with more than 10 GW, 4.65\% exceeding 100 GW, and only 9.30\% requiring less than 1 GW. 

Examples of viability cost curves are presented in Fig. \ref{Fig:boundary_cost_of_4_states}. Essentially, for each considered LDES power capacity, there is an associated amount of total avoided costs that long-duration systems can achieve while supporting thermal replacement. These total avoided costs when divided by the LDES power capacity result in the reported viability costs, below which LDES would be viable (area below the curves). In the graphs of Fig. \ref{Fig:boundary_cost_of_4_states}, initially, for each kW of LDES capacity there is an increasing rate of cost reduction, i.e, higher total avoided costs, which explains growing viability costs at first until a maximum value. As the rate of costs reduction decreases, the viability costs follow the same pattern. 

Figure \ref{Fig:National_capacity_baseline_and_opp_value} shows the total U.S. system capacity and cost by technology in 2050, based on data from the 43 states where the viability costs are positive. Essentially, in case LDES systems are available, gas and coal capacities can be replaced by increasing 43.4\% of solar, 37.7\% of wind-ons, 2\% of wind-ofs, and 0.3\% of SDES capacities, while much more investment would be needed to do the same without LDES systems.     


\subsection{Reasons for high and low viability costs}\label{sec:reasons_high_low_BC}

Different factors and their combinations influence LDES viability costs, leading them to high or low values. Table \ref{Tab:top_10_bottom_10_bc} summarizes the main characteristics of baseline systems (where gas and coal generators are still present and LDES is not included) behind high and low values of LDES viability costs. Essentially, these main characteristics are thermal generation participation in annual load supply, percentage of thermal capacity utilization relative to the available capacity throughout the whole year, average hourly CF of IES in the \textit{baseline model}, gas and coal-related FO\&M costs as a percentage of the overall system costs, and the predominant intermittent generation source (either solar or wind). In addition, we include the parameter $\alpha$, which we define as the ratio between the LDES capacity that maximizes the viability cost and thermal capacity present in the \textit{baseline model}. More information on the results analyzed in this paper can be found in our interactive dashboard \cite{LDES_Tableau} and in the supplementary document.

\begin{table}[h]
\footnotesize
\centering
\caption{Top 10 and bottom 10 viability cost: insights from the \textit{baseline model} run.} \label{Tab:top_10_bottom_10_bc}
    \begin{tabular}{l|l}
        \textbf{Top 10 viability cost} & \textbf{Bottom 10 viability cost} \\ \hline
        \makecell[l]{States: California, Colorado, Idaho, \\ Kansas, Montana, North Dakota, \\ Nebraska, New Mexico, Oklahoma, Virginia} &  
        \makecell[l]{States: Alabama, Connecticut, Delaware, \\ Florida, Kentucky, Maryland, \\ New Jersey, Ohio, Tennessee, Vermont} \\ \hline
        \makecell[l]{9 of the 10 states with the highest \\ 
        viability cost have thermal generation \\  
        participation lower than 10\%} &  
        \makecell[l]{9 of the 10 states with the lowest \\ 
        viability cost have thermal generation \\  
        participation higher than 20\%} \\ \hline
        \makecell[l]{7 of the 10 states with the highest \\
        viability cost have thermal capacity \\ 
        utilization participation lower than 20\%} &  
        \makecell[l]{8 of the 10 states with the lowest \\
        viability cost have thermal capacity \\ 
        utilization participation higher than 20\%} \\ \hline
        \makecell[l]{7 of the 10 states with the highest \\
        viability cost have average hourly capacity \\
        factor (CF) of intermittent energy sources (IES) \\ 
        in baseline higher than 25\%} &  
        \makecell[l]{8 of the 10 states with the lowest \\
        viability cost have average hourly capacity \\ 
        factor (CF) of intermittent energy sources (IES) \\
        in baseline lower than 25\%} \\ \hline
        \makecell[l]{6 of the 10 states with the highest \\ 
        viability cost have thermal FO\&M cost \\
        participation higher than 10\%} &  
        \makecell[l]{7 of the 10 states with the lowest \\
        viability cost have thermal FO\&M cost \\
        participation lower than 10\%} \\ \hline
        \makecell[l]{6 of the 10 states with the highest \\
        viability cost have wind (ons/ofs) generation  \\
        participation in baseline higher than 65\% (\% IES)} &  
        \makecell[l]{8 of the 10 states with the lowest \\
        viability cost have solar generation  \\
        participation in baseline higher than 65\% (\% IES)} \\ \hline
        \makecell[l]{All of the 10 states with the highest \\
        viability cost have $\alpha < 1$} &  
        \makecell[l]{All of the 10 states with the lowest \\
        viability cost have $\alpha > =1$}\\ \hline 
        \end{tabular}
\end{table}

\begin{figure}[!tb]
    \centering
     \includegraphics[width=1.0\textwidth,height=1.0
     \textheight,keepaspectratio]{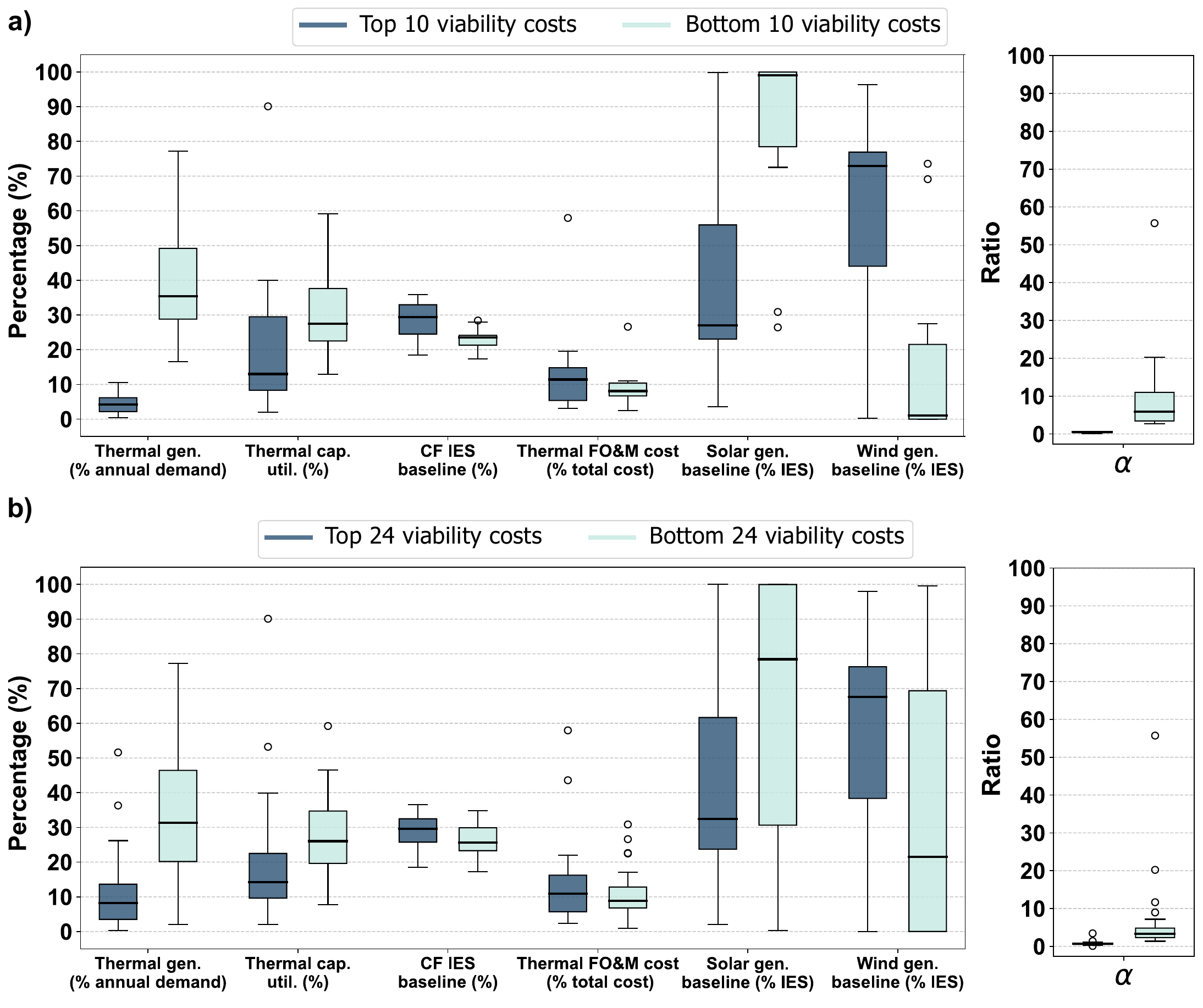}
    \caption{Box plots comparing key system metrics for U.S. states categorized by viability costs. Figure (a) compares the top 10 and bottom 10 states in terms of viability cost, and figure (b) extends to the top 24 and bottom 24, encompassing all states. Metrics include thermal generation participation, thermal capacity utilization, intermittent capacity factors, thermal FO\&M costs, and solar and wind generation in the \textit{baseline model}. The right-side panel presents the distribution of the $\alpha$ ratio for each case.}
    \label{Fig:box_plot_boundary_cost}
\end{figure}



Fig. \ref{Fig:box_plot_boundary_cost} sheds more light on the evaluation of the factors analyzed through thresholds in Table \ref{Tab:top_10_bottom_10_bc} by giving a sense of average and dispersion. Interestingly, even considering top and bottom 24 states in terms of viability costs reaffirms the findings of Table \ref{Tab:top_10_bottom_10_bc}. 
States with lower viability costs usually significantly rely on cheaper thermal generation, whose replacement would require substantial investments in intermittent plants at a comparable or higher cost than maintaining and operating gas and coal units, therefore leaving less room for investments in storage systems. Additionally, since these states normally have a higher proportion of solar capacity in their baseline intermittent mix, they tend to have lower hourly average IES capacity factors in all or almost all meteorological seasons, which also impel more investments in intermittent units to replace the generation provided by thermals and limit the viability of storage systems. On the other hand, higher capacity factors driving periods of excess generation and underutilized thermals with expensive FO\&M costs are associated with higher viability costs.
Concrete examples of this analysis are the cases of Maryland (MD) and North Dakota (ND). MD exhibits a very low LDES viability cost (US\$ 3.94 $\text{kW}^{-1}$), while ND has one of the highest (US\$ 3,881.37 $\text{kW}^{-1}$). MD has a substantial reliance on thermal generation, with 28.37\% of thermal generation participation and a thermal capacity utilization of 22.07\%. 
Additionally, MD's IES exhibits a low CF of 22.05\%, and its IES generation is mostly solar (72.53\%) in the \textit{baseline model}. On the other hand, ND stands at the opposite end of the spectrum, demonstrating the lowest thermal dependence, with a minimal thermal generation share (0.37\%) and an extremely low thermal capacity utilization (6.73\%). Unlike MD, ND's energy mix is dominated by wind, accounting for 73.53\% of its total intermittent generation in the \textit{baseline model}, while solar plays a minor role, 26.47\%.

\subsection{Operation analysis}\label{sec:operation_analysis}

\begin{figure}[!tb]
    \centering
     \includegraphics[width=1.0\textwidth,height=1.0
     \textheight,keepaspectratio]{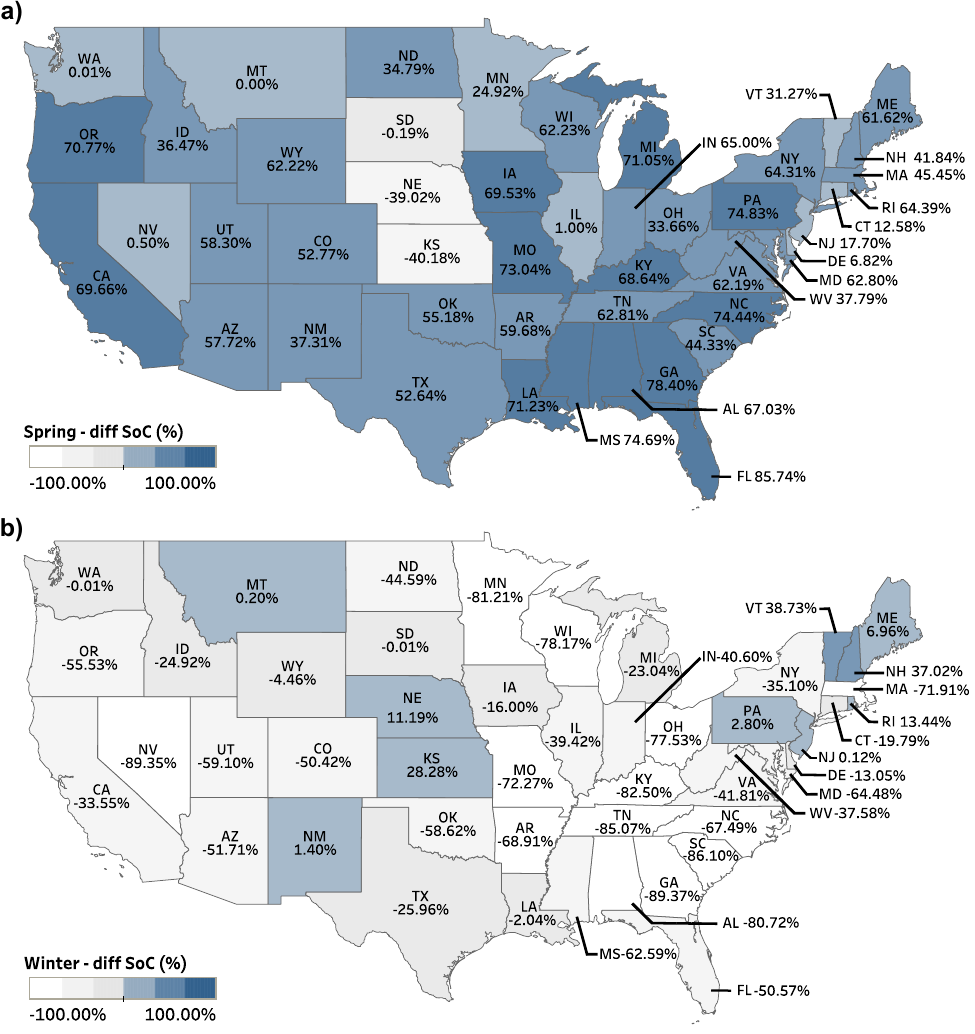}
    \caption{Comparison of seasonal SoC differences: (a) Spring and (b) Winter. During spring, the colors are more intense, indicating higher charging levels, whereas in winter, the colors are lighter, reflecting greater discharging.}
    \label{Fig:soc_dif_seasons}
\end{figure}

\begin{figure}[!tb]
    \centering
    \includegraphics[width=1.0\textwidth,height=1.0
     \textheight,keepaspectratio]
     {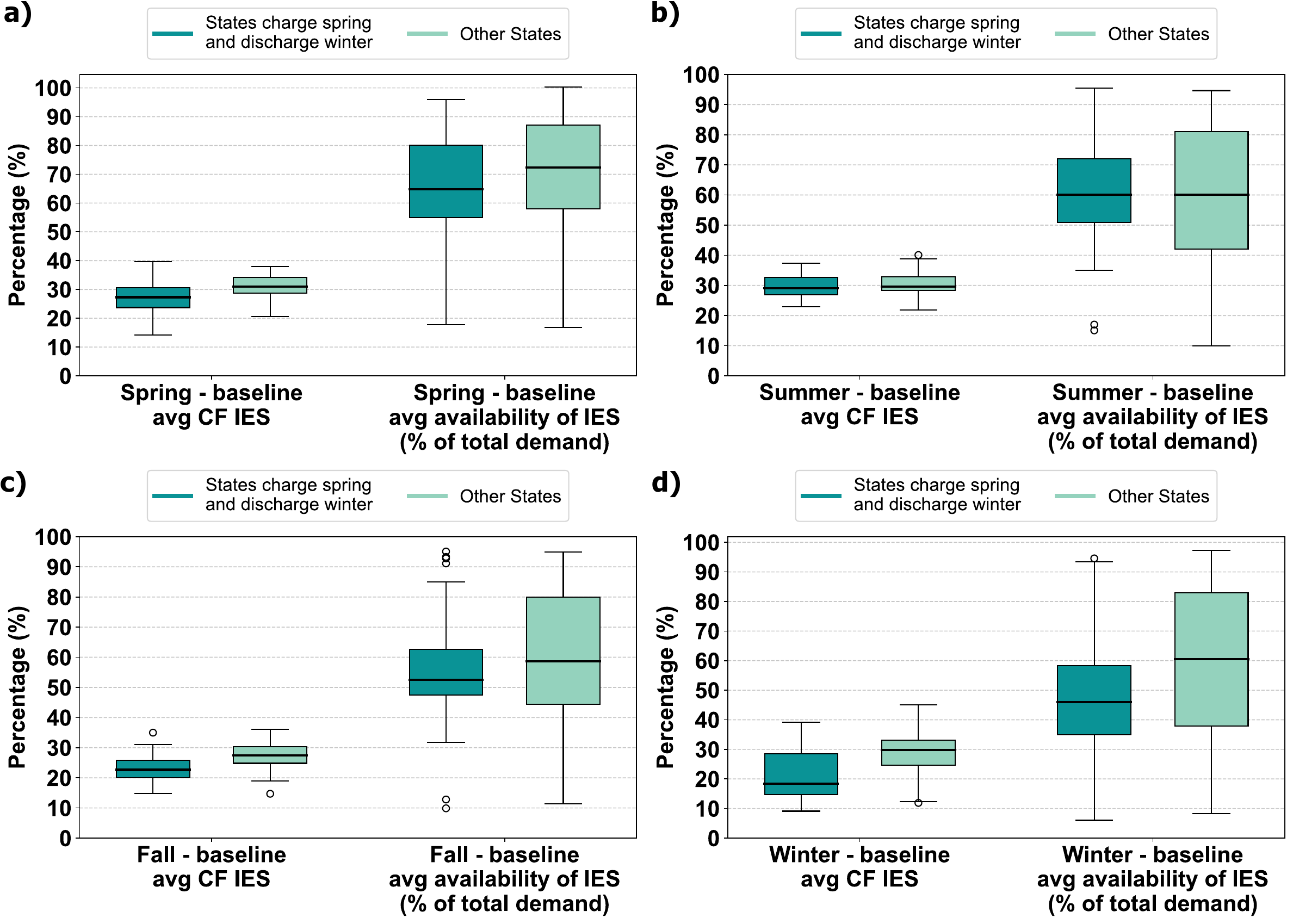}
    \caption{Average capacity factors of IES and average relative availability of IES (in \% of total demand) among the 48 states for each meteorological season.}
    \label{Fig:avg_CF_avg_rel_IES_avail_per_season}
\end{figure}

\begin{figure}[!tb]
    \centering
    \includegraphics[width=1.0\textwidth,height=1.0
     \textheight,keepaspectratio]
     {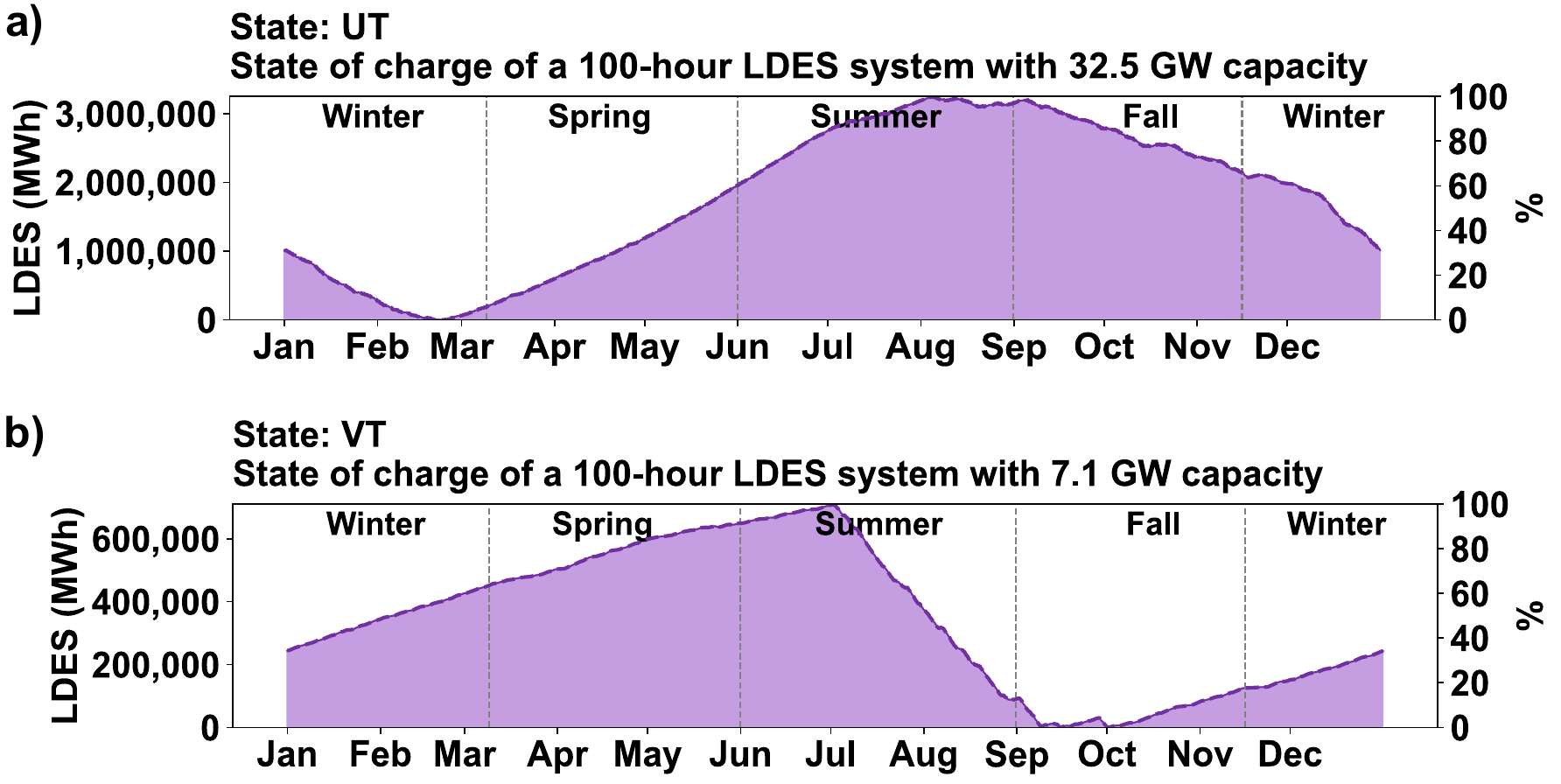}
    \caption{Different LDES SOC patterns. (a) LDES SOC through the year in Utah. (b) LDES SOC through the year in Vermont.}
    \label{Fig:soc_2_states_UT_and_VT}
\end{figure}

The impact of LDES systems varies across states depending on the season, as depicted in Fig. \ref{Fig:soc_dif_seasons} for spring and winter (other seasons illustrated in the supplementary document). This figure presents the seasonal relative differences in State-of-Charge (SoC), computed as the SoC at the end of a season minus its initial SoC normalized by the SoC capacity for each state. 

Our analysis indicates that 36 states would charge more during spring, while 6 in the summer, 4 in the fall, and 2 in the winter. Moreover, 28 states would discharge more in the winter, 12 in the summer, 5 in the fall, and 3 in the spring. In fact, 22 states charge more in the spring and discharge more in the winter. This pattern occurs since average relative availability of IES (\% of total demand) is particularly higher in spring and lower during winter for these 22 states, while other states in general also have differences but less exacerbated as shown by Fig. \ref{Fig:avg_CF_avg_rel_IES_avail_per_season}. Two contrasting examples are Utah and Vermont, whose annual state of charge profiles are illustrated in Fig. \ref{Fig:soc_2_states_UT_and_VT}. Utah belongs to the aforementioned group of 22 states and therefore leverages its higher availability of intermittent plants to charge in the spring aiming to discharge in the winter when IES generation is less available. On the other hand, Vermont exhibits an unusual behavior charging more during winter making use of the especially higher wind speeds in part of the Northeastern US during this season \cite{Klink1999}.

\section{Conclusion}\label{sec:conclusion}

We investigated in this paper the state-wise economic viability of LDES from the system perspective in each state given its characteristics. Several factors impact this viability. In general, states with higher LDES viability costs tend to have a lower participation of thermal units for supplying the annual demand along with under-utilization of these resources, and a higher share of wind plants in their intermittent portfolio combined with higher CFs. Given the future projection of energy matrices and LDES cost targets, we found that the economic viability of 100-h LDES systems to replace conventional thermal plants will be limited to 4 states.  

\section{Methods}\label{sec:method}

\subsection{Methodology}\label{sec:sensitivity_analysis}

Our methodology to conduct the studies in this paper leverages available public data and the optimization models formulated in \cite{first_paper}, with the objective to determine the viability cost of LDES, i.e., the cost below this technology becomes economically viable in the 48 states of the contiguous US. The optimization models are a \textit{baseline model} and an \textit{opportunity value maximization model} which are solved in sequence. Firstly, the \textit{baseline model} is solved with the goal of minimizing both operational and investment expenditures associated with the selected system. Following this, given the cost discovered, the \textit{opportunity value maximization model} is solved and aims to maximize LDES's opportunity value while keeping the same cost of \textit{baseline model}. In the \textit{opportunity value maximization model}, an investment portfolio of intermittent and SDES assets is optimized to replace existing user-selected generation technologies, given the inclusion of a user-defined LDES power capacity (we assume 100 hours duration in our studies). Each imposed LDES power capacity might be able to avoid a particular amount of costs associated with investments and operational decisions to retire a particular set of technologies (e.g., gas and coal). We define the maximum \textit{opportunity value} for LDES as the maximum avoided costs that can be achieved. The maximum LDES viability cost (in \$/kW) is then defined as the ratio between the maximum \textit{opportunity value} and the LDES capacity associated with it.

\subsection{\textit{Baseline model}}\label{sec:baseline_model}

The \textit{baseline bodel} aims to minimize the total system cost, as described in Eq. \eqref{eq:baseline_objective}. The cost components include investment costs for generation and short-duration storage (\(C^{\text{inv,gen}}\) and \(C^{\text{inv,st,short}}\)), operational costs (\(C^{\text{ope,gen}}\)), reserve costs (\(C^{\text{ies,gen}}\)), imbalance costs (\(C^{\text{imbalance}}\)), reserve shortage costs (\(C^{\text{ies,shortage}}\)), and fixed operation and maintenance costs for generation and storage (\(C^{\text{FO\&M,gen}}\) and \(C^{\text{FO\&M,st,short}}\)). The optimization is subject to a set of constraints: power balance, reserve margin requirements, operational limits of storage devices, and technical limits of generators. This model establishes a reference system cost under ideal operational conditions.

\begin{align}
q^* = \min \Big( 
    & \hspace{5pt} C^{\text{inv,gen}} + C^{\text{inv,st,short}} \notag \\
    & \hspace{5pt} + C^{\text{ope,gen}} + C^{\text{ies,gen}} + C^{\text{imbalance}} + C^{\text{ies,shortage}} \notag \\
    & \hspace{5pt} + C^{\text{FO\&M,gen}} + C^{\text{FO\&M,st,short}} 
\Big)
\label{eq:baseline_objective}
\end{align}

\subsection{\textit{Opportunity value maximization model}}\label{sec:opp_value_max_model}

The aim of the \textit{opportunity value maximization model} is described in Eq. \eqref{eq:opportunity_objective}. It maximizes the difference between the variable called viability cost of LDES (\(c^{\text{VC}}\)) and the penalization (\(C^{\text{over}} q^{\text{over}}\)) of potentially incurring in a total system cost higher than that of the baseline model. This model is subject to the same set of constraints of the \textit{baseline model} in addition to an extra constraint \eqref{eq:const_opp_value_model} that limits the total cost of operation plus investment, including the new LDES (\(c^{\text{VC}} x^{\text{power}}\)) capacity, to be lower or equal to the reference cost $q^*$ obtained from the baseline model. In case this limit cannot be satisfied, the relaxing auxiliary variable $q^{\text{over}}$, associated with over cost, assumes value greater than zero, which is the case for the states with negative viability costs. 

\begin{equation}
q^* = \max \Big( c^{\text{VC}} - C^{\text{over}} q^{\text{over}} \Big)
\label{eq:opportunity_objective}
\end{equation}

\begin{equation}
\begin{aligned}
\text{s.t.:} \hspace{5pt}
& \Big( C^{\text{inv,gen}} + C^{\text{inv,st,short}} \\
& \hspace{5pt} + C^{\text{ope,gen}} + C^{\text{ies,gen}} + C^{\text{imbalance}} + C^{\text{ies,shortage}} \\
& \hspace{5pt} + C^{\text{FO\&M,gen}} + C^{\text{FO\&M,st,short}} \\
& \hspace{5pt} + c^{\text{VC}} x^{\text{power}} \Big) 
\leq q^* + q^{\text{over}}
\end{aligned}
\label{eq:const_opp_value_model}
\end{equation}

\subsection{Assumptions in the models for this study}

In the \textit{baseline model}, although investment variables are present (essentially to have the same set of general decision variables of the \textit{opportunity value maximization model}), specific constraints are included to prevent investments, i.e., only the existing baseline resources (in our case study, the data for these resources come from NREL's Cambium projection for 2050 \cite{gagnon2023cambium_site}) can be dispatched and no investments are allowed in LDES, SDES, or any type of new generator. Additionally, there is no retirement of existing firm generators, ensuring that the system operates with its given infrastructure. In contrast, in the \textit{opportunity value maximization model}, LDES capacity is incorporated based on a predetermined value provided by the user. Investments in SDES and intermittent energy sources are also permitted, and the retirement of gas and coal generators is considered.

\subsection{Data source}

The dataset combines projections of the USA's 2050 state-wise energy matrix. System load and generators, including installed capacities, are based on NREL's Cambium 2022 dataset (no tax credit phaseout scenario) \cite{gagnon2023cambium_site}. Specific characteristics of technologies, such as ramp rates, come from Cambium 2022 documentation \cite{gagnon2023cambium}, while fixed O\&M costs are sourced from the ReEDS base \cite{2021reeds} and the 2022 Annual Technology Baseline (ATB) \cite{atb2022}. Investment costs for intermittent energy are also derived from the ATB 2022. Fuel prices are reported in AEO 2023 \cite{AEO_2023}. CFs of intermittent generation were calculated using data from Cambium 2022 and ReEDS base \cite{2021reeds}. Our financial results in this study are expressed in 2022 U.S. dollars.

\subsection{Spatial and temporal resolution}

The simulation model was conducted for the entire United States (48 states), with each state analyzed individually. Power plants of different technologies are located within specific balancing areas (BAs), with 134 BAs across the U.S. \cite{gagnon2023cambium}. The analysis for each state uses an hourly temporal resolution, totaling 8,760 intervals over the year, where the combined generating capacity of all BAs within a state must meet the state's demand. To ensure system reliability, we impose a reserve requirement of 4\% of the system demand across all time periods and states. This value aligns with standard practices in grid planning and operation, where reserve margins are maintained to accommodate fluctuations in generation capacity availability and demand.

\subsection{Fixed and candidate generators}

The database provided by Cambium 2022 includes projections of existing generators across the United States for the year 2050. In this study, K-means clustering was employed to reduce the number of generators by grouping them into three representative categories, high, medium, and low cost, within each BA, thereby decreasing computational complexity. Apart from the existing generators, we also considered a set of candidate intermittent generators. For each intermittent technology (solar, wind-ons, and wind-ofs) within a BA, a corresponding candidate generator with identical specifications was introduced. The maximum investment in each candidate generator was set to 4x (four times) the installed capacity of the corresponding fixed generator, ensuring flexibility for system expansion. However, for specific states like Connecticut (CT), Delaware (DE), and Pennsylvania (PA), this multiplier was increased to 10x (ten times), reflecting these states’ unique demand profiles and geographic limitations, which necessitate greater capacity expansion to achieve a feasible solution.

\subsection{Fixed and candidate batteries}

In this study, fixed batteries in each BA were grouped by type with durations of 2h, 4h, 6h, 8h, and pumped hydro storage (PHS). For candidate short-duration energy storage (SDES), a single 4-hour system was created for each state, with a maximum investment capacity set to 10 times the total existing storage capacity of SDES technologies in the state (excluding PHS). Additionally, a candidate LDES system, also for each state, with a 100-hour duration was included, with capacity limits ranging from 0 to 150 GW and an round trip efficiency of 42.5\%. The SDES round trip efficiency was set at 85\%.

\input{References.bbl}



\section*{Supplementary material}

# 1 Supplementary table with key parameters analyzed in the paper

Key parameter values analyzed in the paper for all states, including the power quantity that maximizes the viability cost, are available in Table 1.

Table 1: All values of key parameters for U.S. states.

| State | Maximum viability cost ($/kW) | Power qty. (GW) for max viability cost | Thermal gen. (% annual demand) | Thermal cap. util. (%) | CF IES baseline (%) | Thermal FO&M cost (% total cost) | Solar gen. baseline (% total IES) | Wind gen. baseline (% total IES) | α |
|---|---|---|---|---|---|---|---|---|---|
| OH | -11.81 | 150.00 | 50.56 | 59.21 | 23.09 | 6.65 | 96.32 | 3.68 | 4.55 |
| NJ | -3.95 | 150.00 | 30.15 | 28.89 | 28.41 | 7.44 | 26.43 | 73.57 | 11.66 |
| CT | -1.43 | 150.00 | 56.31 | 38.97 | 17.85 | 8.94 | 98.07 | 1.93 | 20.22 |
| AL | -0.60 | 150.00 | 31.97 | 33.60 | 23.88 | 8.75 | 100.00 | 0.00 | 8.96 |
| DE | -0.07 | 150.00 | 77.27 | 23.82 | 17.30 | 10.86 | 100.00 | 0.00 | 55.71 |
| MD | 3.94 | 36.40 | 28.37 | 22.07 | 22.05 | 11.01 | 72.53 | 27.47 | 3.73 |
| KY | 22.73 | 75.00 | 45.07 | 15.88 | 23.86 | 26.61 | 99.96 | 0.04 | 3.34 |
| TN | 23.60 | 39.80 | 16.55 | 26.15 | 24.23 | 2.52 | 99.90 | 0.10 | 3.22 |
| FL | 27.17 | 116.00 | 38.90 | 46.57 | 21.07 | 6.85 | 100.00 | 0.00 | 2.67 |
| VT | 50.83 | 7.10 | 20.51 | 12.90 | 27.90 | 5.14 | 30.87 | 69.13 | 7.17 |
| UT | 97.72 | 32.50 | 30.78 | 26.63 | 29.50 | 30.85 | 98.50 | 1.50 | 5.65 |
| AR | 101.39 | 30.20 | 13.89 | 7.69 | 23.36 | 22.49 | 100.00 | 0.00 | 2.41 |
| IA | 104.43 | 19.70 | 12.18 | 32.58 | 31.46 | 8.42 | 26.90 | 73.10 | 2.95 |
| MA | 114.28 | 29.90 | 36.49 | 39.13 | 30.45 | 9.93 | 29.87 | 70.13 | 4.52 |
| MI | 118.40 | 79.30 | 53.70 | 41.85 | 31.61 | 7.67 | 33.95 | 66.05 | 2.63 |
| PA | 119.46 | 130.00 | 53.89 | 38.02 | 26.45 | 11.46 | 29.89 | 70.11 | 3.31 |
| SC | 129.59 | 37.90 | 19.15 | 21.98 | 20.10 | 5.27 | 100.00 | 0.00 | 2.07 |
| LA | 134.57 | 36.40 | 42.65 | 19.01 | 24.87 | 16.97 | 84.34 | 15.66 | 1.63 |
| IN | 139.73 | 38.50 | 53.50 | 22.62 | 29.69 | 17.03 | 41.61 | 58.39 | 1.36 |
| WY | 148.94 | 13.20 | 3.01 | 13.78 | 34.87 | 10.43 | 0.37 | 99.63 | 3.96 |
| GA | 149.06 | 45.00 | 20.91 | 25.92 | 24.91 | 5.87 | 100.00 | 0.00 | 2.01 |
| IL | 165.57 | 6.00 | 1.99 | 19.86 | 32.95 | 0.97 | 31.28 | 68.72 | 1.47 |
| NY | 173.75 | 32.90 | 24.24 | 29.79 | 27.10 | 7.21 | 49.89 | 50.11 | 1.45 |
| RI | 178.96 | 10.30 | 41.19 | 18.38 | 32.17 | 22.62 | 27.91 | 72.09 | 3.01 |
| WI | 179.89 | 17.90 | 36.31 | 11.85 | 25.18 | 21.95 | 66.63 | 33.37 | 0.80 |
| AZ | 217.58 | 12.90 | 11.81 | 16.14 | 23.23 | 12.54 | 98.26 | 1.74 | 0.86 |
| MS | 227.35 | 18.60 | 26.18 | 23.17 | 28.78 | 15.61 | 53.40 | 46.60 | 1.39 |
| MO | 257.84 | 14.00 | 7.65 | 7.44 | 36.50 | 19.80 | 13.72 | 86.28 | 0.73 |
| SD | 260.82 | 0.70 | 0.59 | 14.18 | 33.67 | 2.31 | 2.04 | 97.96 | 0.65 |
| OR | 266.02 | 3.00 | 3.46 | 22.28 | 32.33 | 6.84 | 39.44 | 60.56 | 1.07 |
| NH | 288.79 | 4.40 | 21.64 | 14.22 | 32.20 | 7.08 | 25.62 | 74.38 | 0.98 |
| NV | 309.70 | 17.40 | 23.69 | 53.22 | 30.31 | 5.94 | 100.00 | 0.00 | 3.43 |
| WV | 328.92 | 10.00 | 51.60 | 4.92 | 26.77 | 43.59 | 31.90 | 68.10 | 0.69 |
| MN | 333.95 | 3.00 | 8.80 | 24.78 | 28.36 | 4.69 | 40.65 | 59.35 | 0.53 |
| WA | 353.88 | 8.00 | 1.14 | 4.99 | 30.72 | 10.38 | 15.34 | 84.66 | 0.64 |
| NC | 384.30 | 10.00 | 13.56 | 18.77 | 24.58 | 5.04 | 77.98 | 22.02 | 0.47 |
| ME | 393.54 | 1.30 | 10.16 | 13.84 | 31.46 | 18.12 | 32.87 | 67.13 | 0.78 |
| TX | 445.05 | 46.50 | 13.90 | 19.64 | 25.95 | 11.06 | 32.05 | 67.95 | 0.61 |
| OK | 480.30 | 5.40 | 3.52 | 14.23 | 32.96 | 12.19 | 14.95 | 85.05 | 0.52 |
| NM | 659.73 | 1.20 | 5.98 | 33.83 | 31.05 | 3.11 | 27.67 | 72.33 | 0.56 |
| CO | 743.27 | 2.80 | 1.72 | 10.27 | 34.55 | 14.17 | 22.75 | 77.25 | 0.52 |
| ID | 801.35 | 3.70 | 4.62 | 2.04 | 18.48 | 57.96 | 99.80 | 0.20 | 0.56 |
| CA | 859.88 | 7.70 | 6.23 | 16.40 | 32.93 | 6.95 | 59.99 | 40.01 | 0.28 |
| MT | 977.11 | 0.30 | 0.37 | 11.77 | 35.89 | 4.85 | 3.61 | 96.39 | 0.66 |
| KS | 1512.43 | 2.40 | 10.60 | 39.93 | 26.13 | 10.66 | 43.82 | 56.18 | 0.46 |
| NE | 3793.01 | 0.50 | 3.69 | 7.67 | 27.68 | 19.51 | 24.05 | 75.95 | 0.10 |
| ND | 3881.37 | 0.50 | 0.37 | 6.73 | 19.20 | 15.01 | 26.47 | 73.53 | 0.22 |
| VA | 5993.94 | 1.00 | 9.09 | 90.11 | 23.97 | 3.73 | 79.25 | 20.75 | 0.67 |

4

# 2 Supplementary figure on seasonal SoC charging/discharging

Figure 1 specifically highlights the comparison of seasonal SoC differences for Summer and Fall.

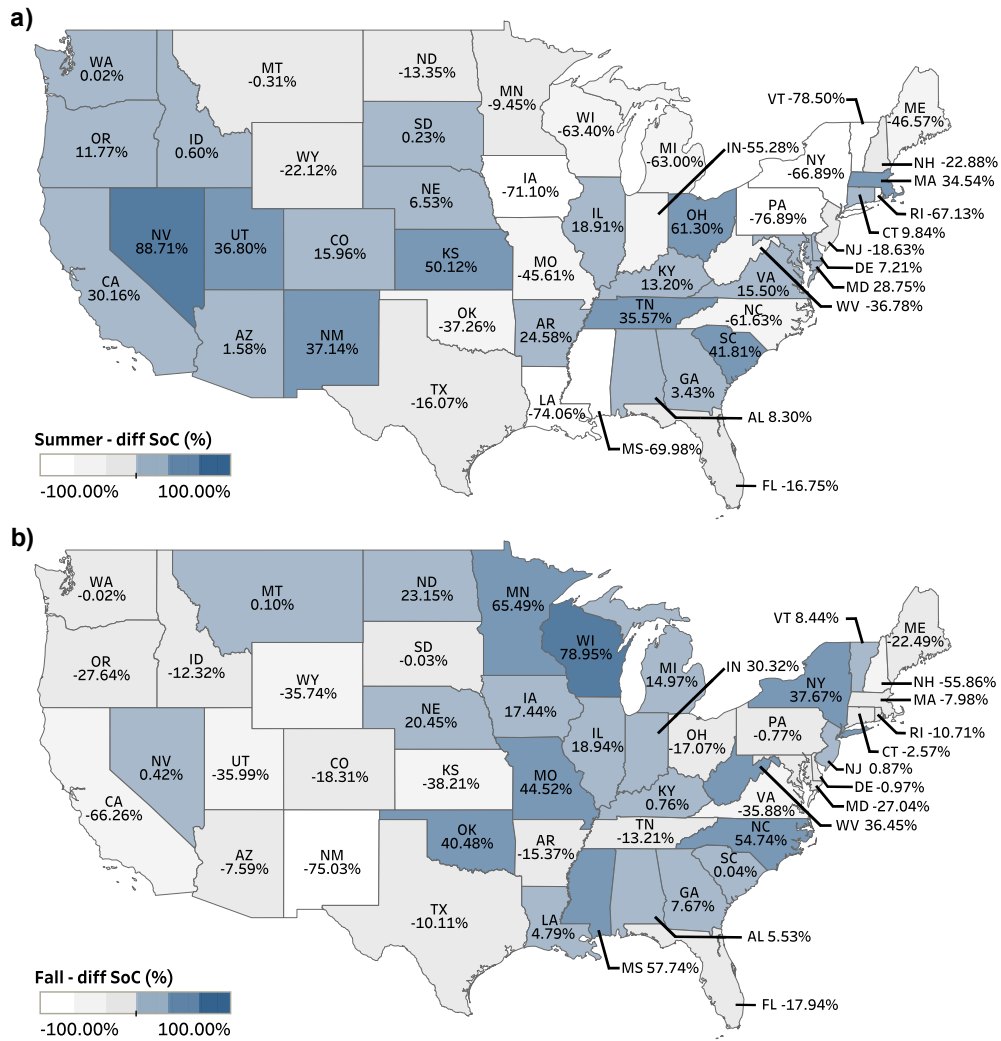

**Fig. 1**: Comparison of seasonal SoC differences: (a) Summer and (b) Fall.

# 3 Supplementary graphs on the operation of all the states of the United States.

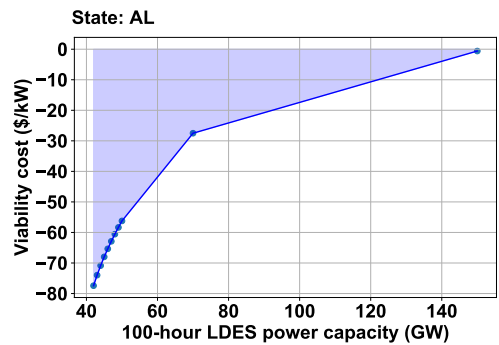
(a) Viability cost curve.

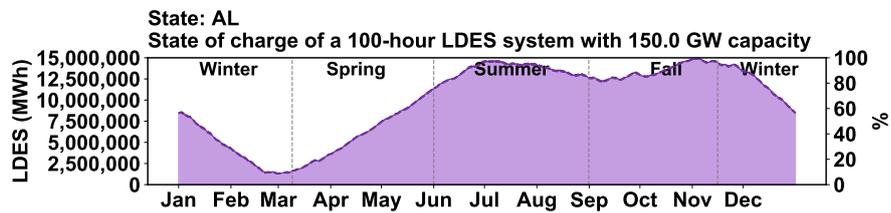
(b) State of charge.

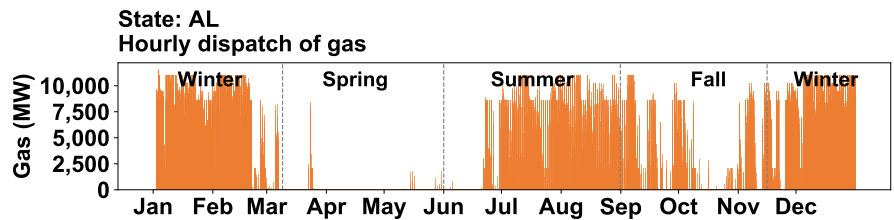
(c) Hourly dispatch of gas.

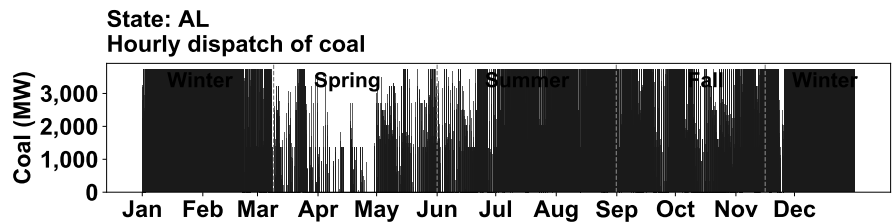
(d) Hourly dispatch of coal.

**Fig. 2**: Alabama (AL).

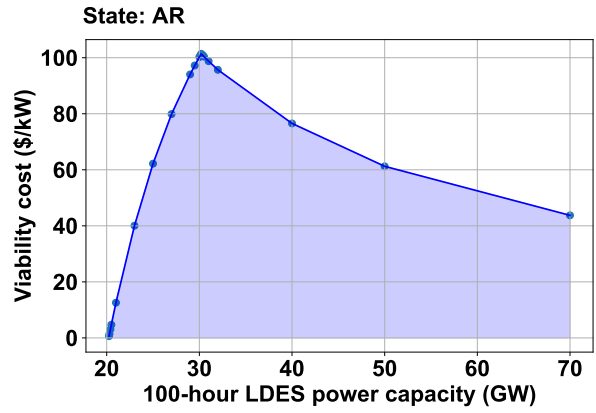

(a) Viability cost curve.

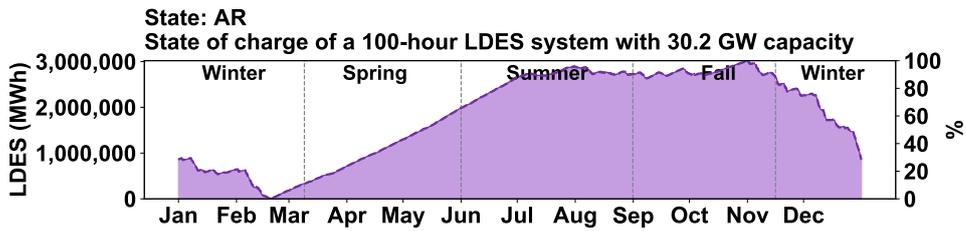

(b) State of charge.

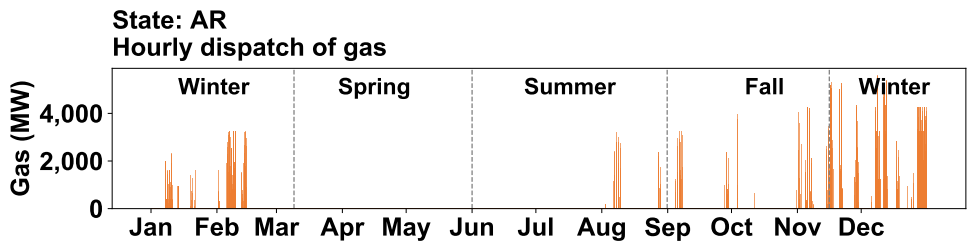

(c) Hourly dispatch of gas.

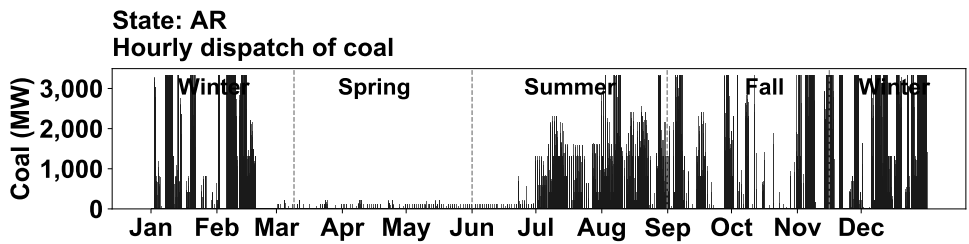

(d) Hourly dispatch of coal.

**Fig. 3**: Arkansas (AR).

7

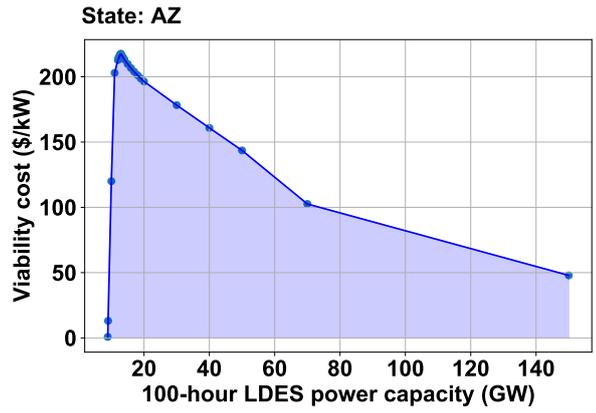

(a) Viability cost curve.

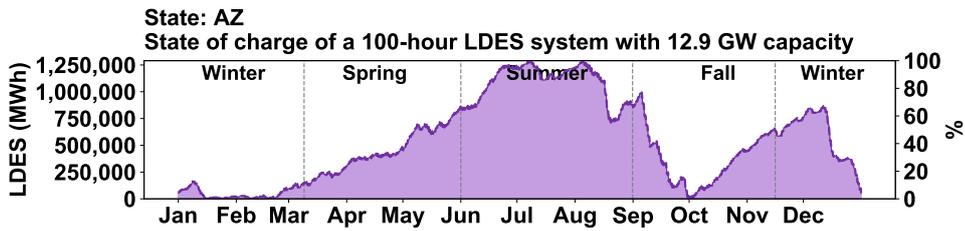

(b) State of charge.

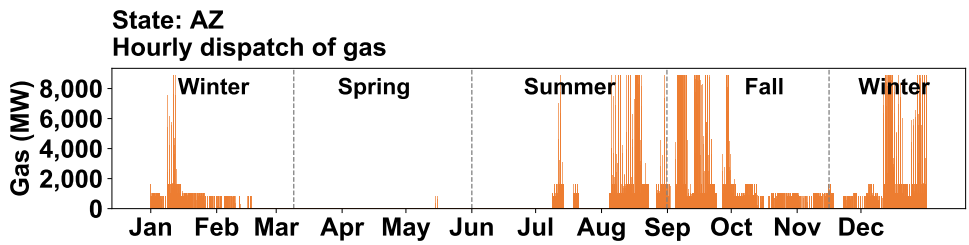

(c) Hourly dispatch of gas.

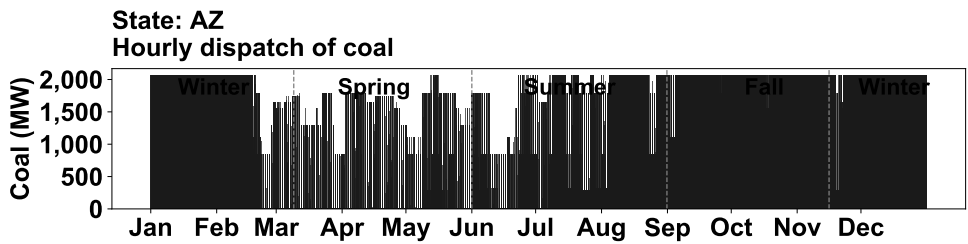

(d) Hourly dispatch of coal.

**Fig. 4**: Arizona (AZ).

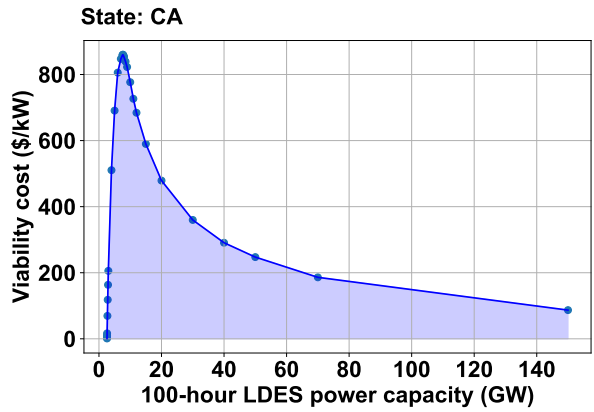

(a) Viability cost curve.

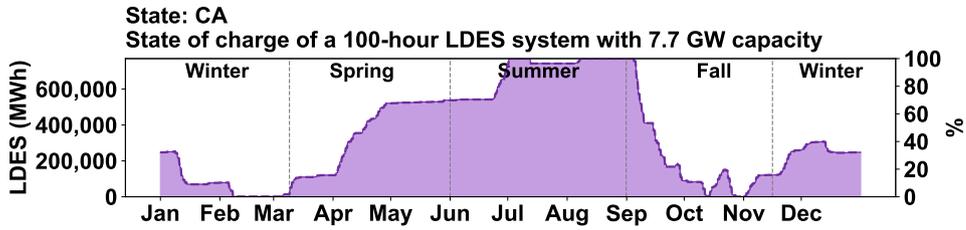

(b) State of charge.

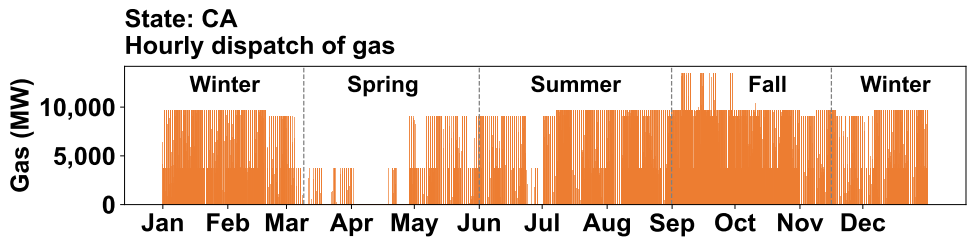

(c) Hourly dispatch of gas.

**Fig. 5**: California (CA).

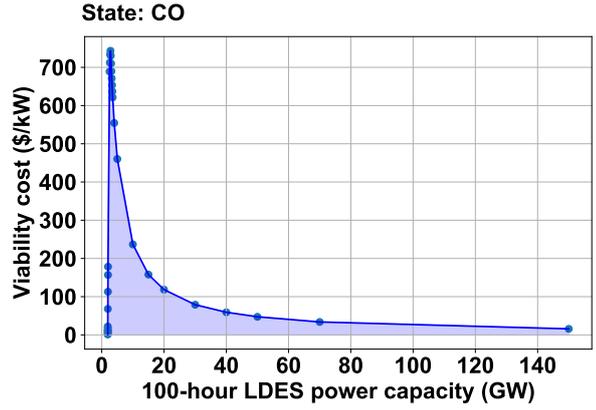

(a) Viability cost curve.

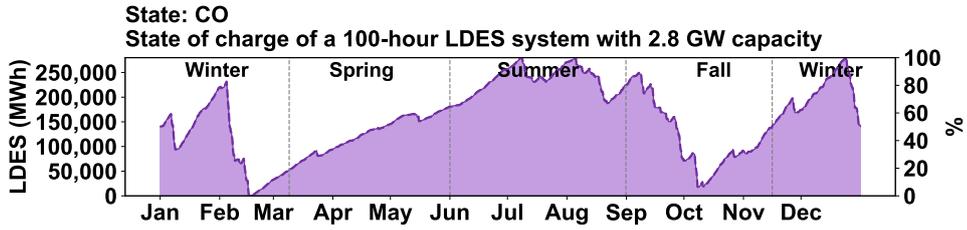

(b) State of charge.

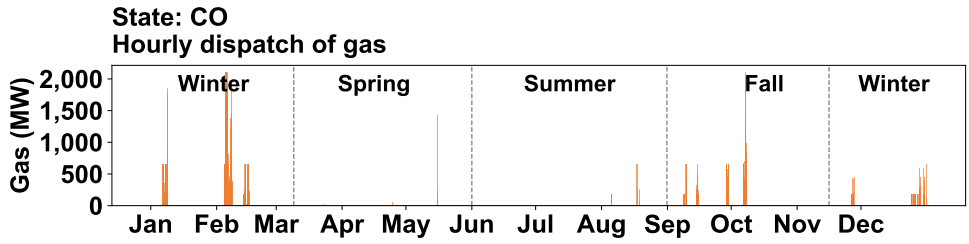

(c) Hourly dispatch of gas.

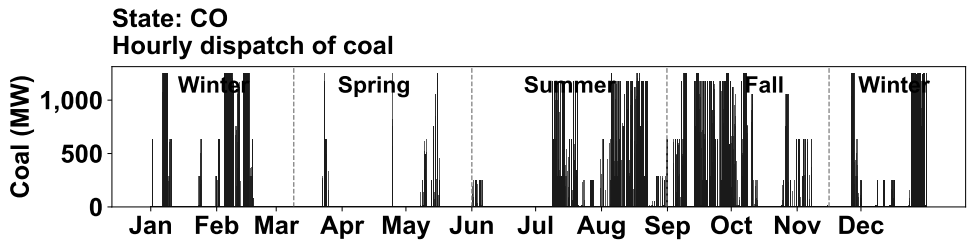

(d) Hourly dispatch of coal.

**Fig. 6**: Colorado (CO).

10

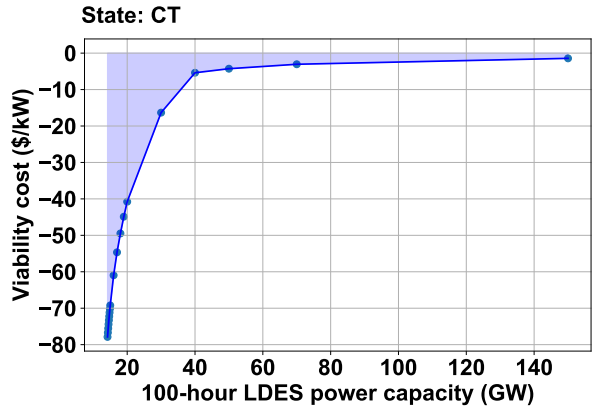

(a) Viability cost curve.

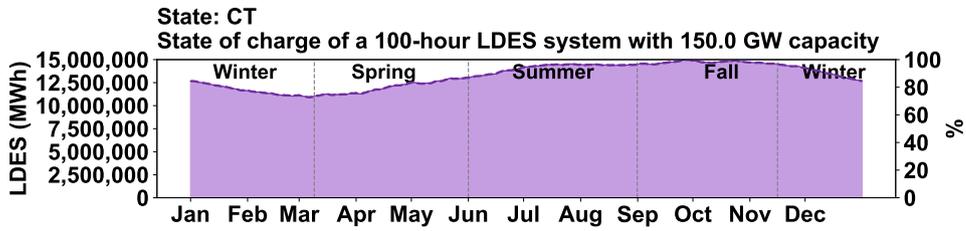

(b) State of charge.

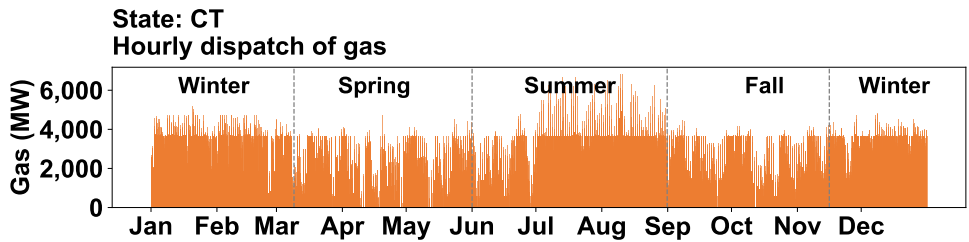

(c) Hourly dispatch of gas.

**Fig. 7**: Connecticut (CT).

11

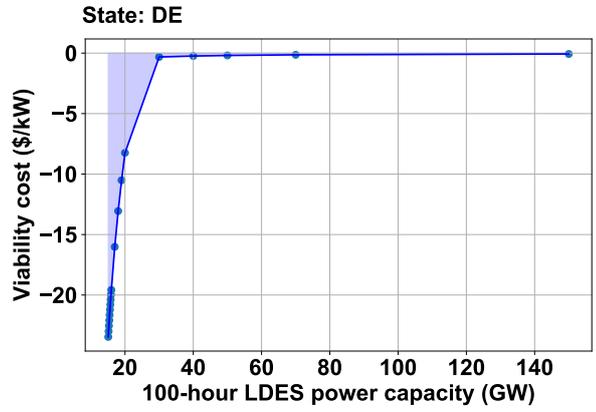

(a) Viability cost curve.

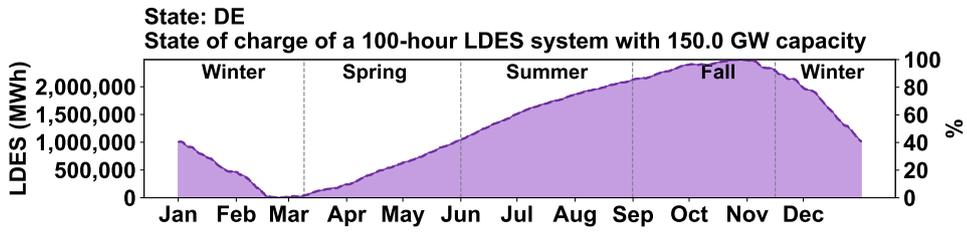

(b) State of charge.

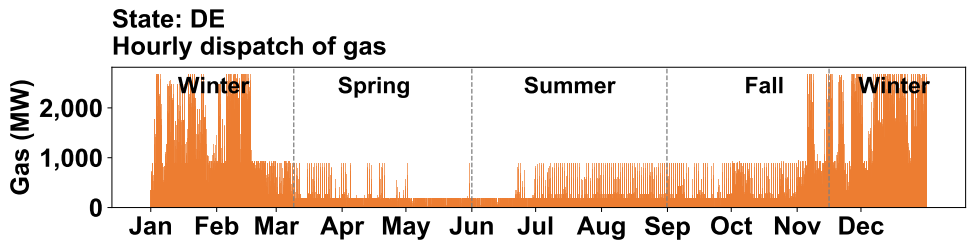

(c) Hourly dispatch of gas.

**Fig. 8**: Delaware (DE).

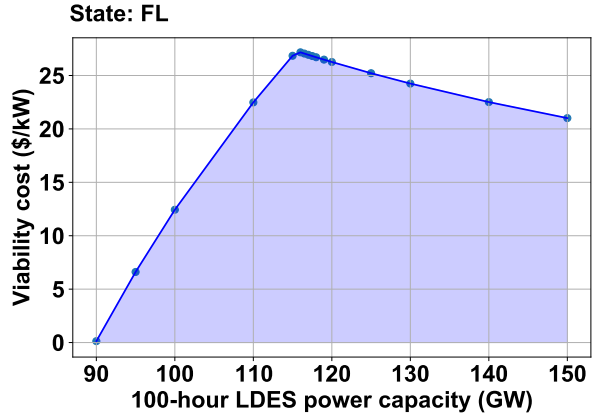

(a) Viability cost curve.

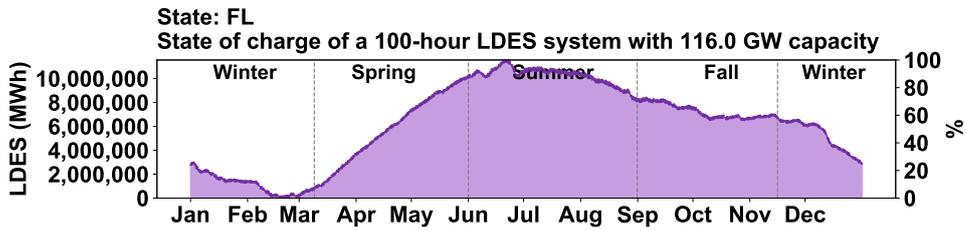

(b) State of charge.

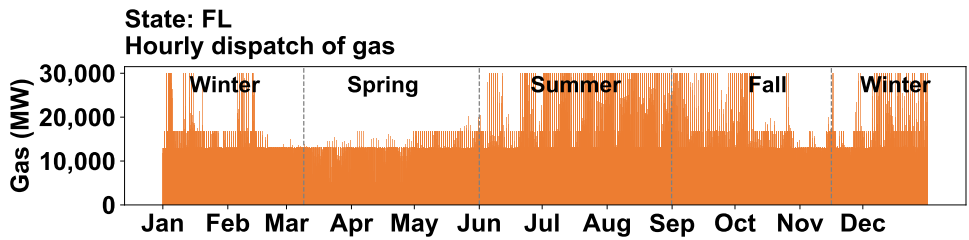

(c) Hourly dispatch of gas.

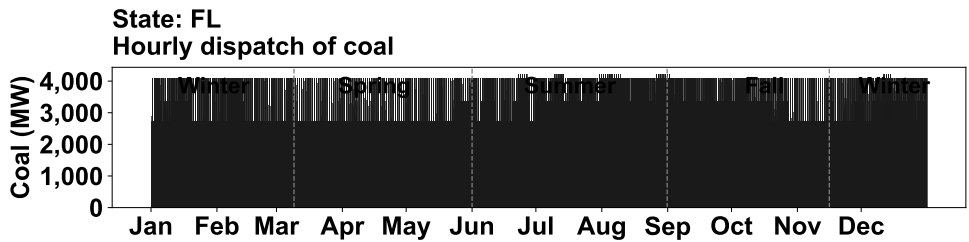

(d) Hourly dispatch of coal.

**Fig. 9**: Florida (FL).

13

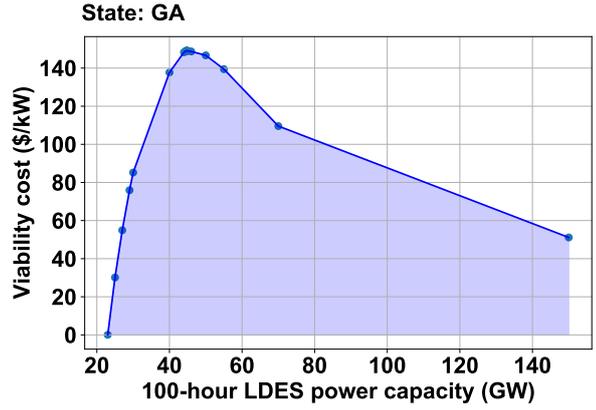

(a) Viability cost curve.

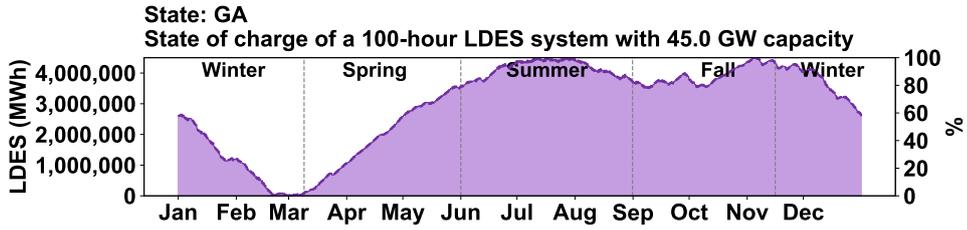

(b) State of charge.

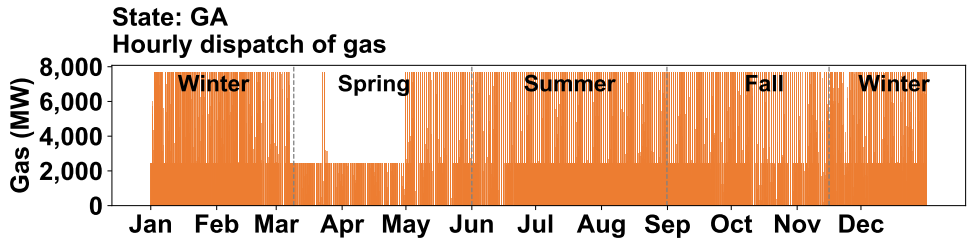

(c) Hourly dispatch of gas.

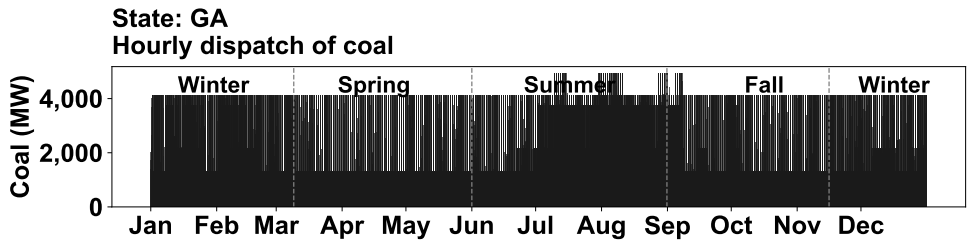

(d) Hourly dispatch of coal.

**Fig. 10**: Georgia (GA).

14

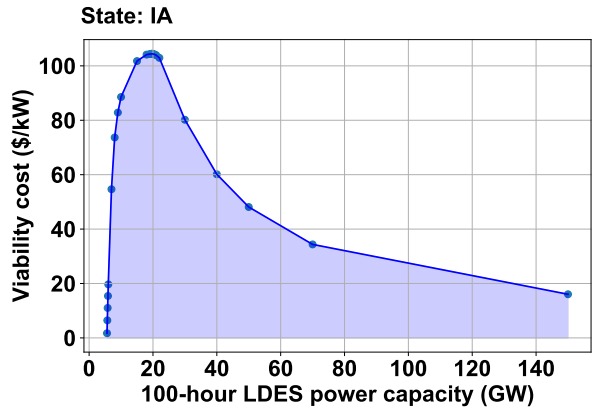

(a) Viability cost curve.

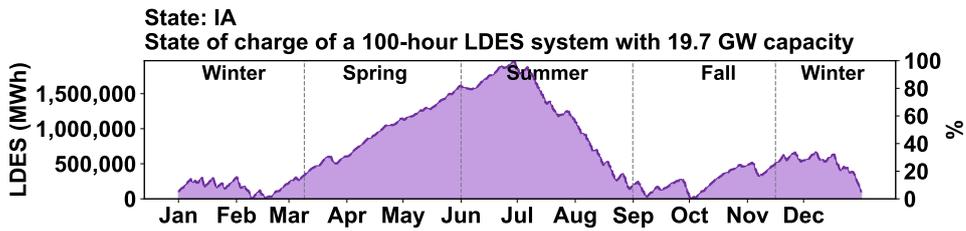

(b) State of charge.

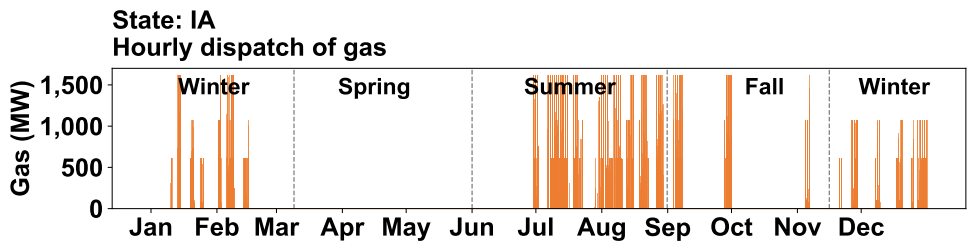

(c) Hourly dispatch of gas.

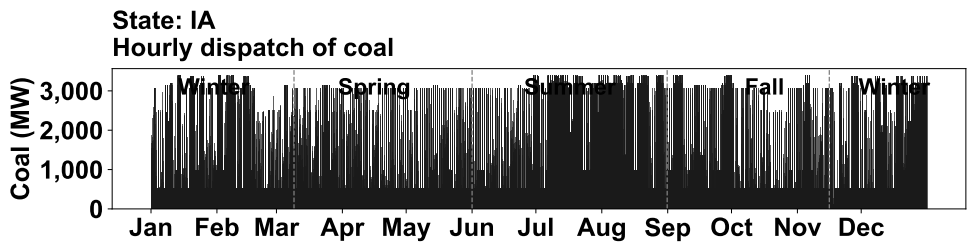

(d) Hourly dispatch of coal.

**Fig. 11**: Iowa (IA).

15

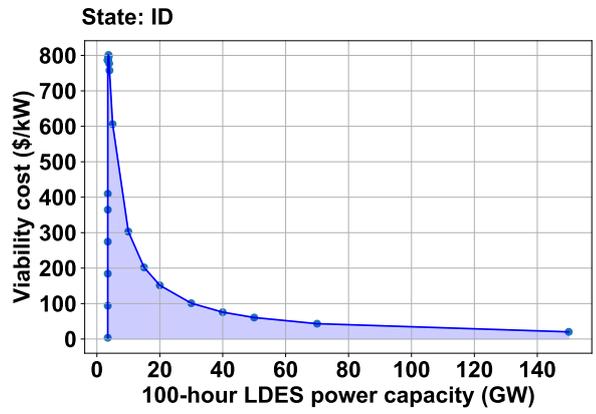

(a) Viability cost curve.

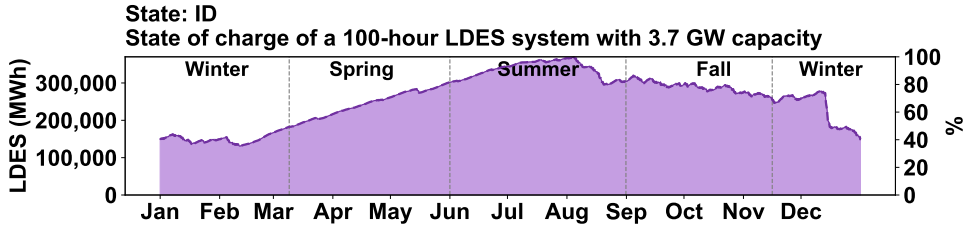

(b) State of charge.

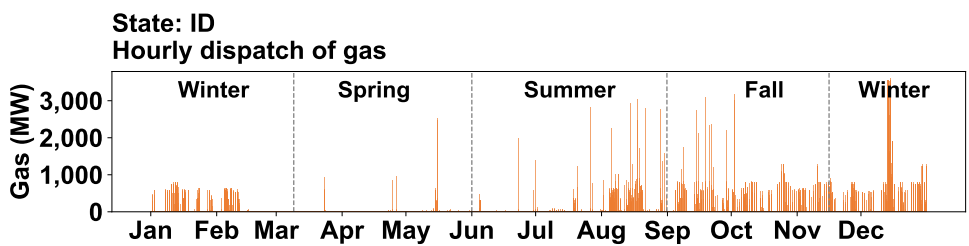

(c) Hourly dispatch of gas.

**Fig. 12**: Idaho (ID).

16

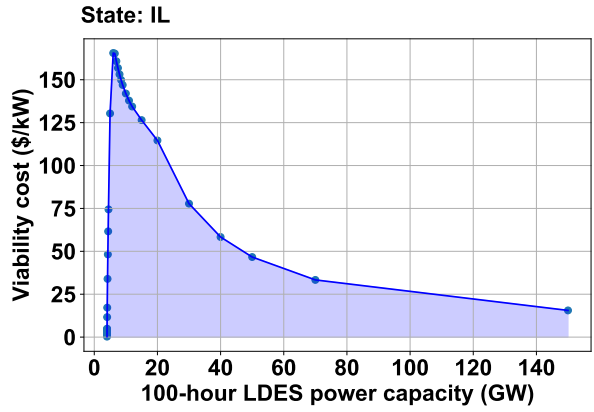

(a) Viability cost curve.

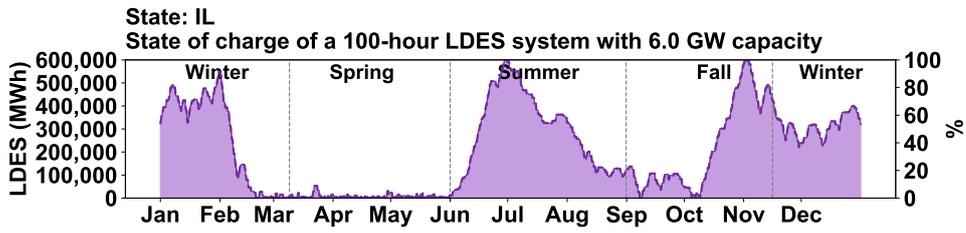

(b) State of charge.

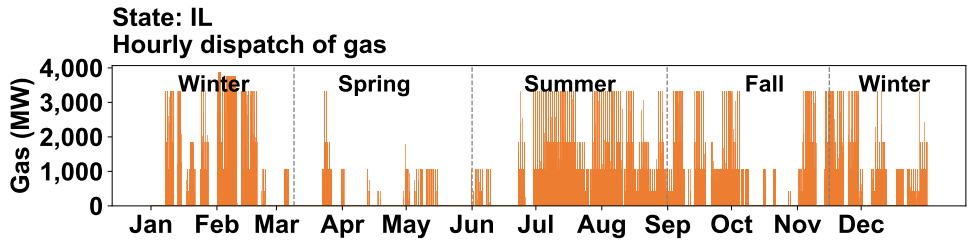

(c) Hourly dispatch of gas.

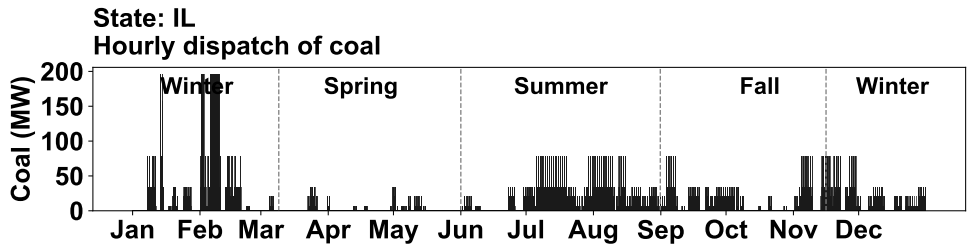

(d) Hourly dispatch of coal.

**Fig. 13**: Illinois (IL).

17

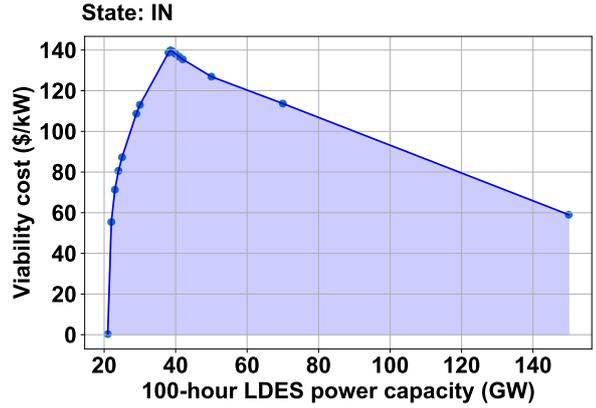

(a) Viability cost curve.

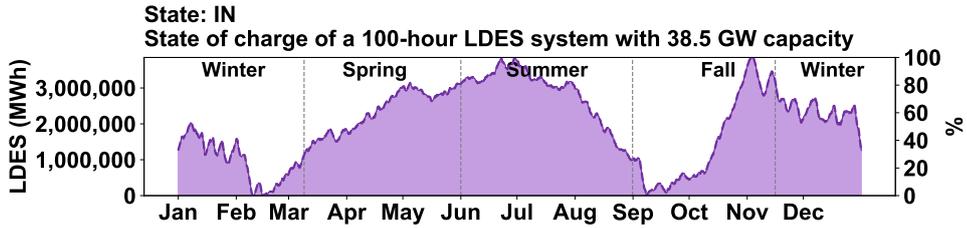

(b) State of charge.

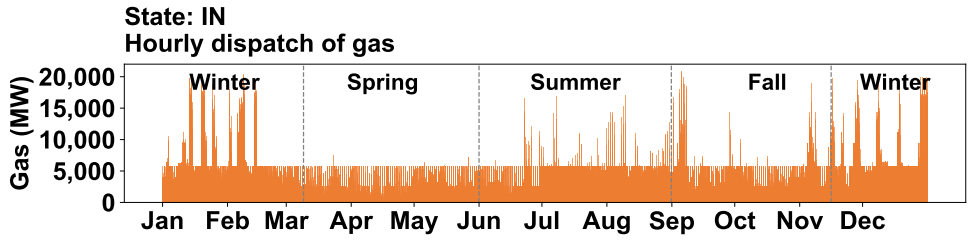

(c) Hourly dispatch of gas.

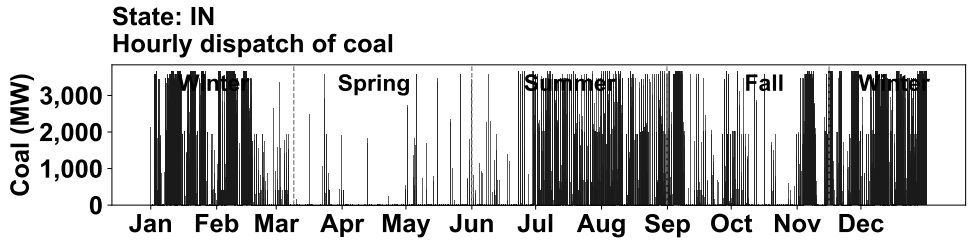

(d) Hourly dispatch of coal.

**Fig. 14**: Indiana (IN).

18

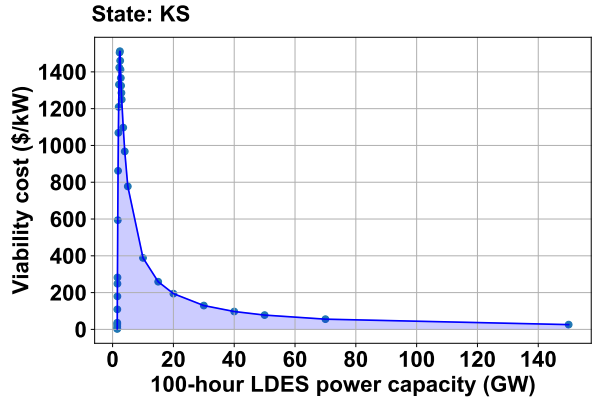
(a) Viability cost curve.

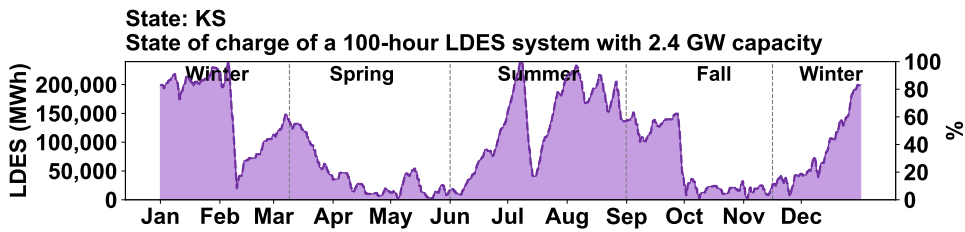
(b) State of charge.

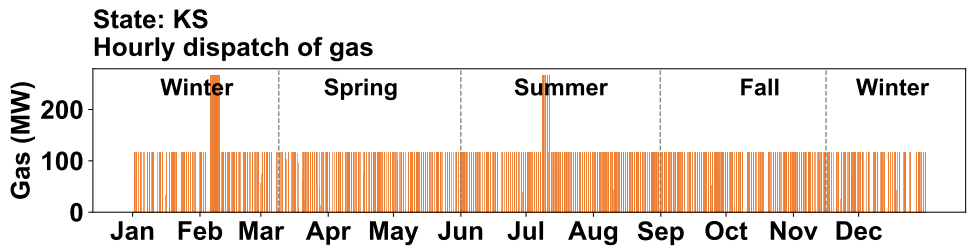
(c) Hourly dispatch of gas.

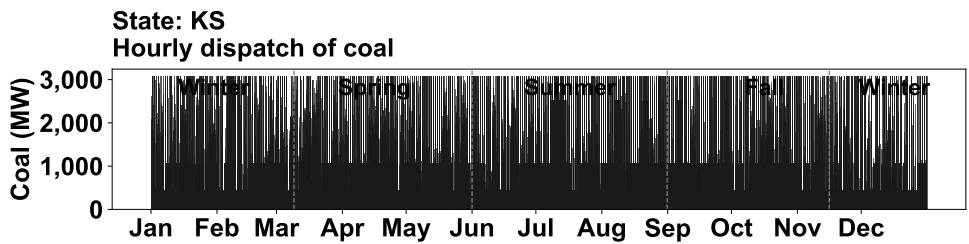
(d) Hourly dispatch of coal.

**Fig. 15**: Kansas (KS).

19

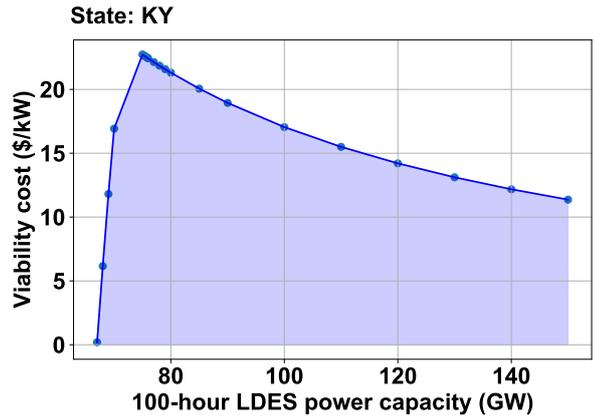

(a) Viability cost curve.

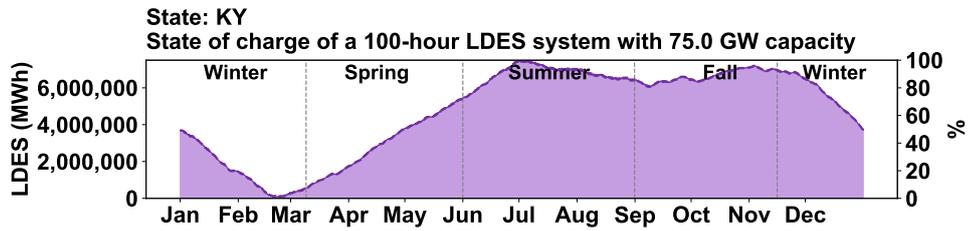

(b) State of charge.

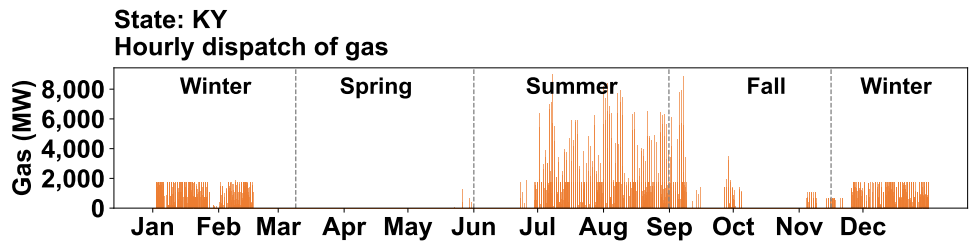

(c) Hourly dispatch of gas.

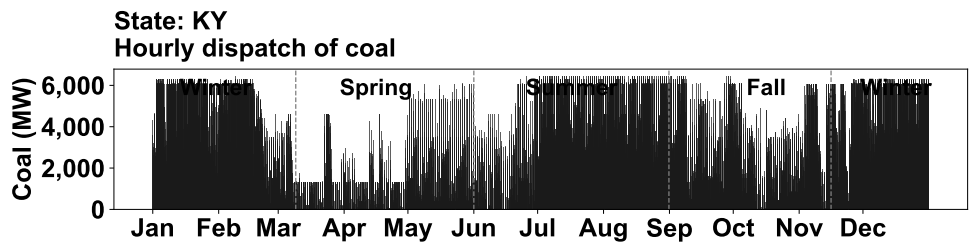

(d) Hourly dispatch of coal.

**Fig. 16**: Kentucky (KY).

20

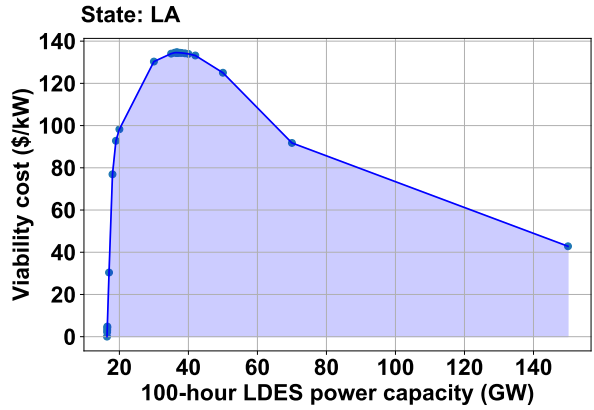

(a) Viability cost curve.

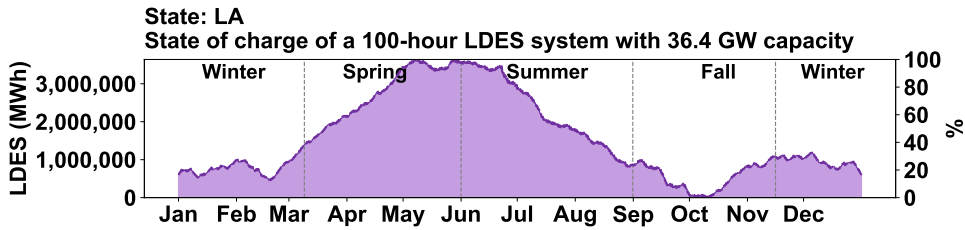

(b) State of charge.

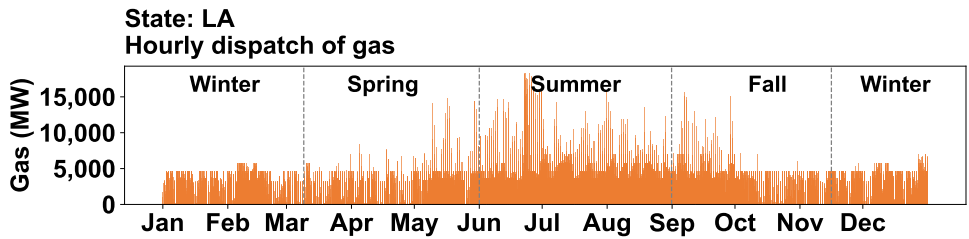

(c) Hourly dispatch of gas.

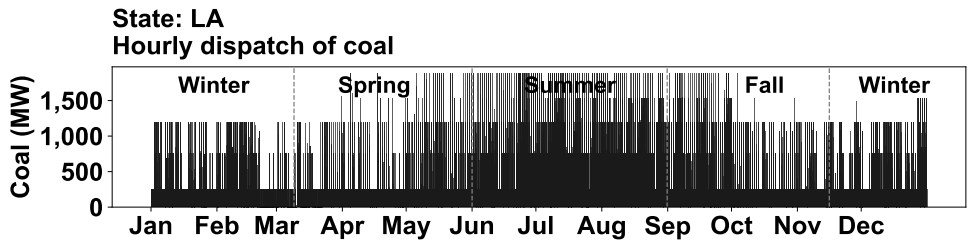

(d) Hourly dispatch of coal.

**Fig. 17**: Louisiana (LA).

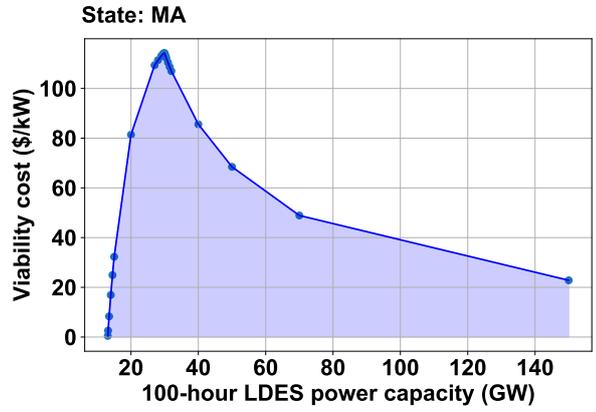

(a) Viability cost curve.

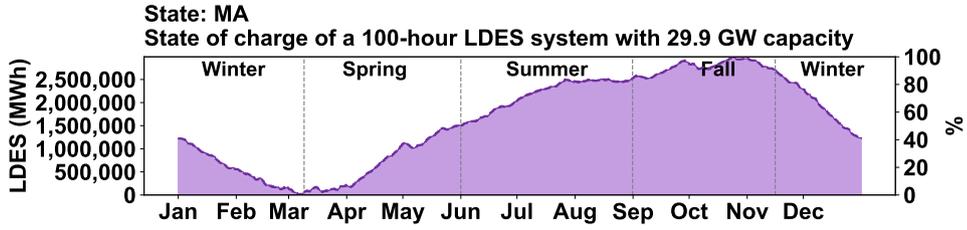

(b) State of charge.

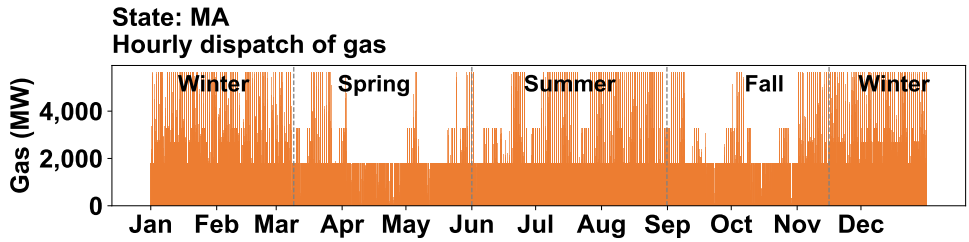

(c) Hourly dispatch of gas.

**Fig. 18**: Massachusetts (MA).

22

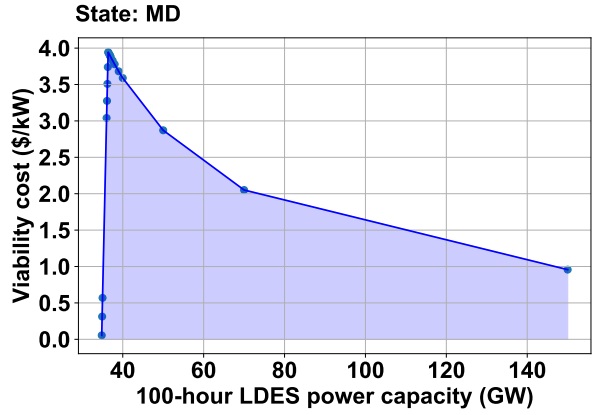

(a) Viability cost curve.

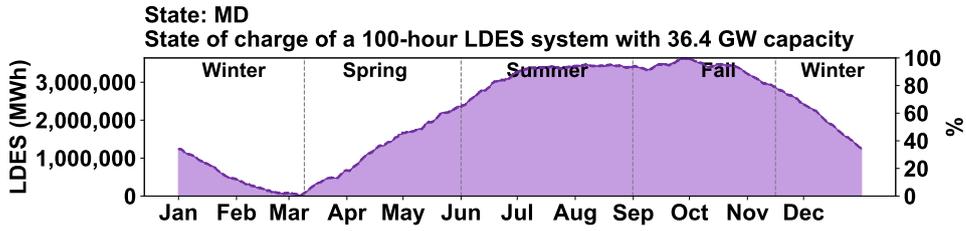

(b) State of charge.

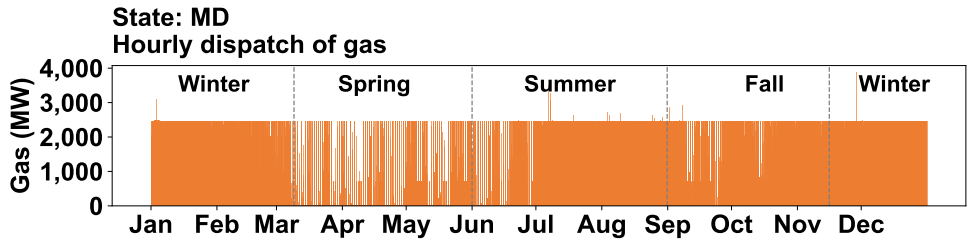

(c) Hourly dispatch of gas.

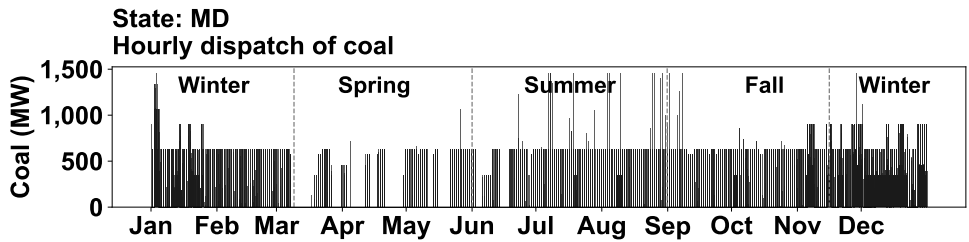

(d) Hourly dispatch of coal.

**Fig. 19**: Maryland (MD).

23

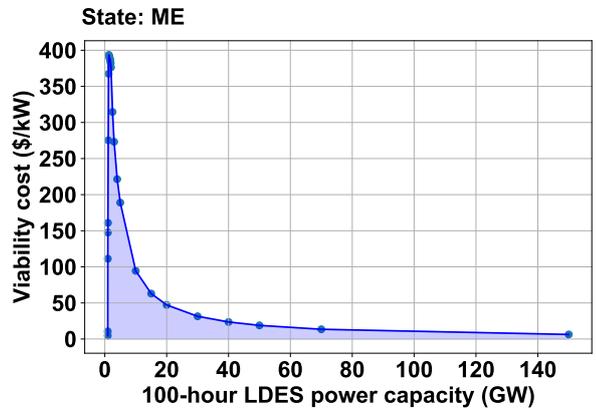

(a) Viability cost curve.

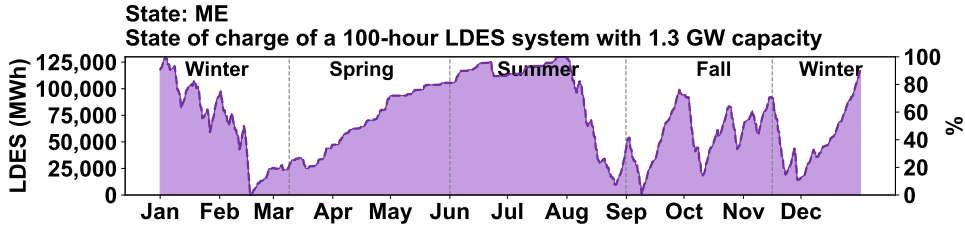

(b) State of charge.

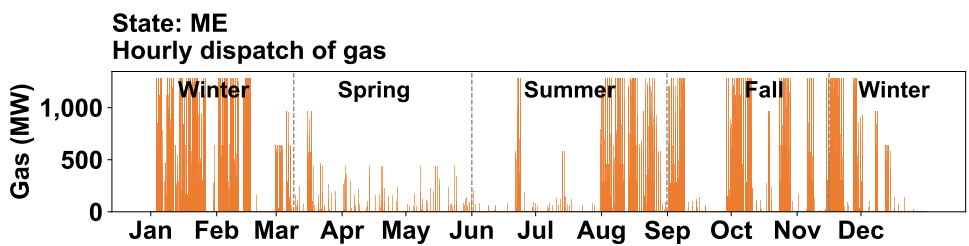

(c) Hourly dispatch of gas.

**Fig. 20**: Maine (ME).

24

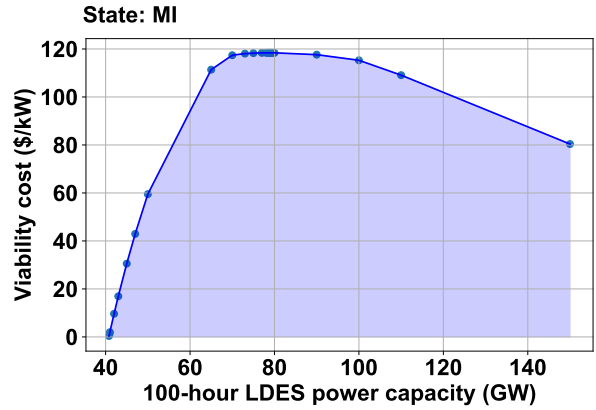

(a) Viability cost curve.

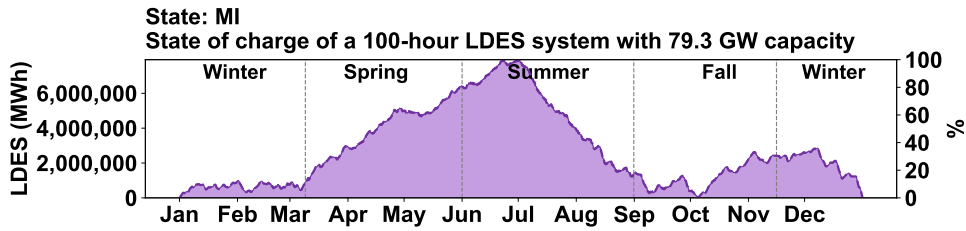

(b) State of charge.

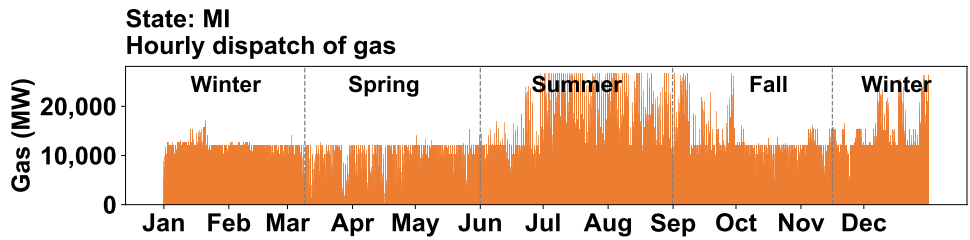

(c) Hourly dispatch of gas.

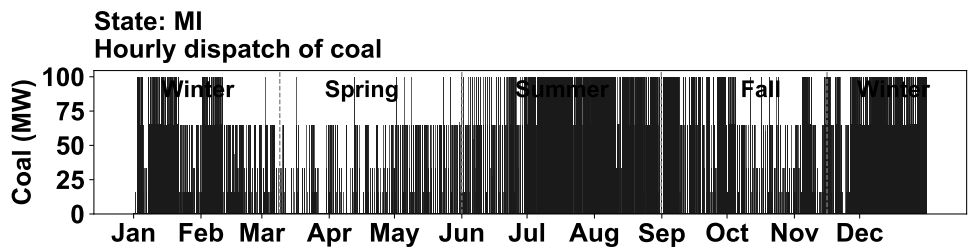

(d) Hourly dispatch of coal.

**Fig. 21**: Michigan (MI).

25

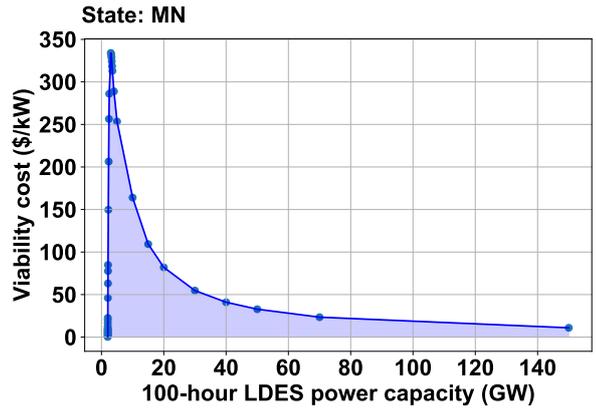

(a) Viability cost curve.

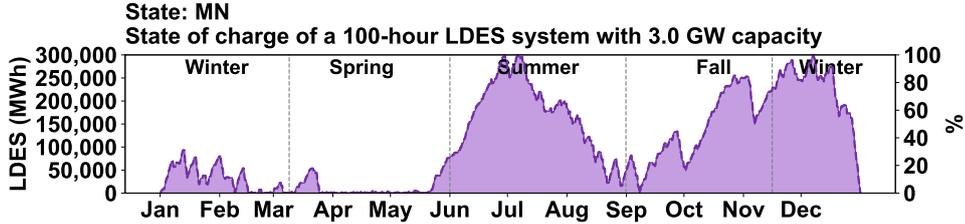

(b) State of charge.

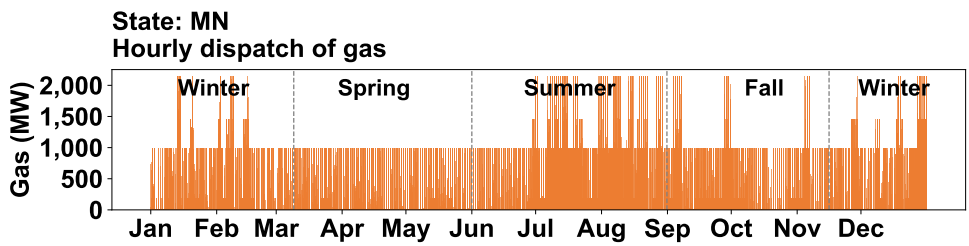

(c) Hourly dispatch of gas.

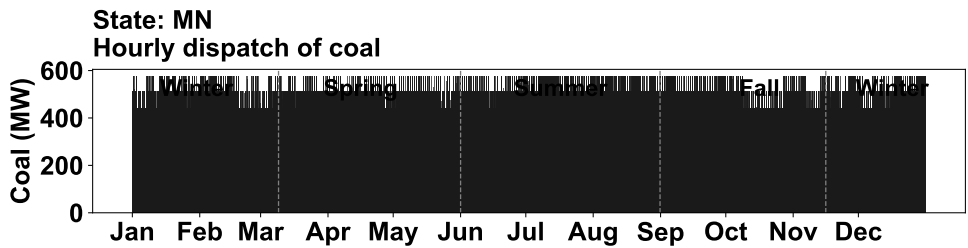

(d) Hourly dispatch of coal.

**Fig. 22**: Minnesota (MN).

26

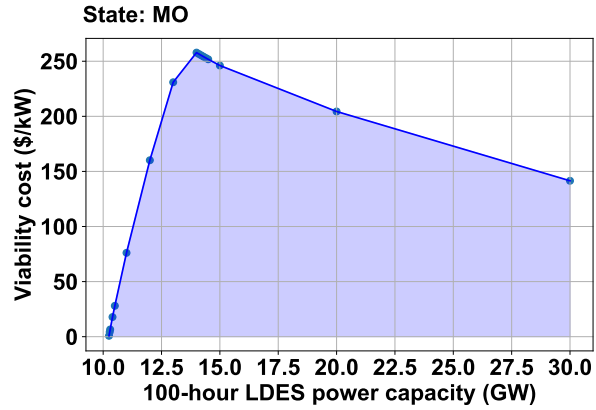
(a) Viability cost curve.

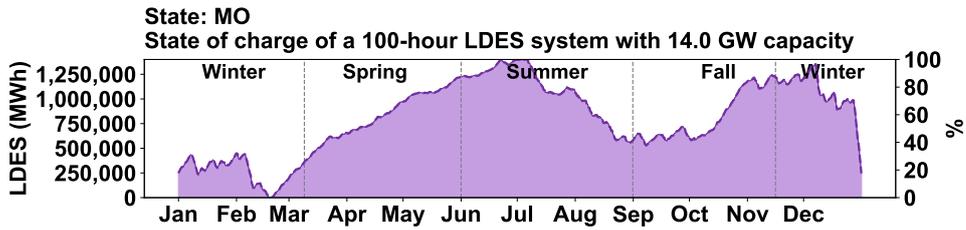
(b) State of charge.

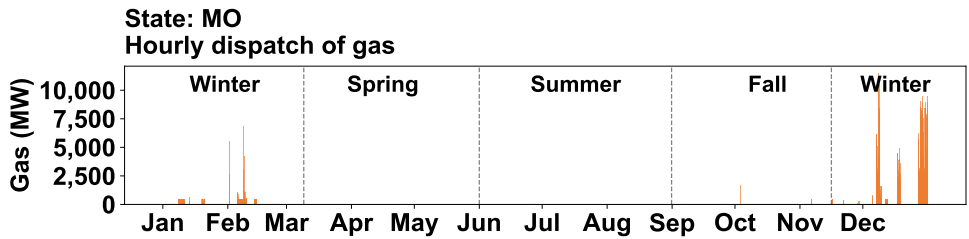
(c) Hourly dispatch of gas.

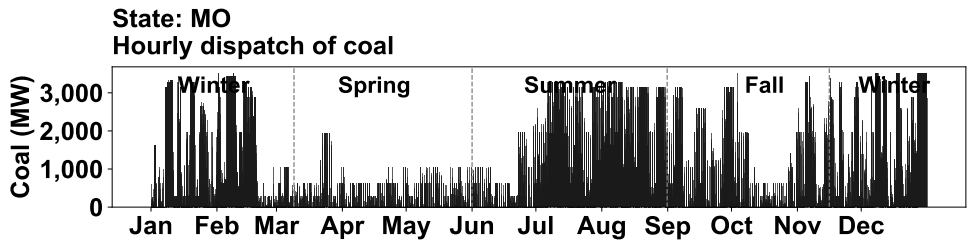
(d) Hourly dispatch of coal.

Fig. 23: Missouri (MO).

27

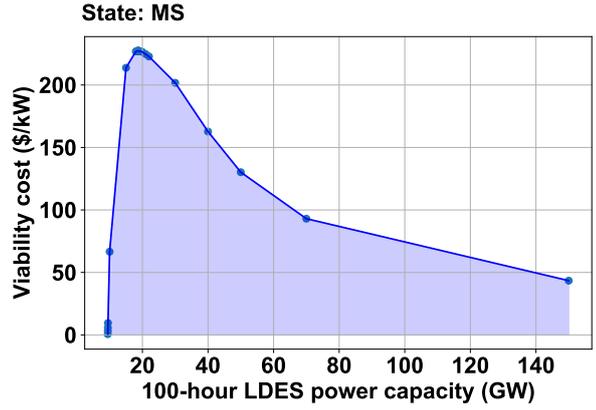

(a) Viability cost curve.

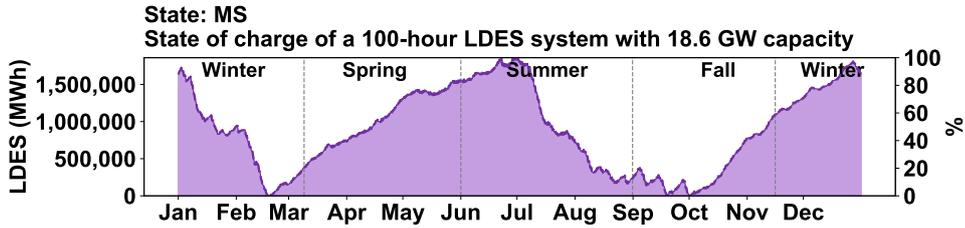

(b) State of charge.

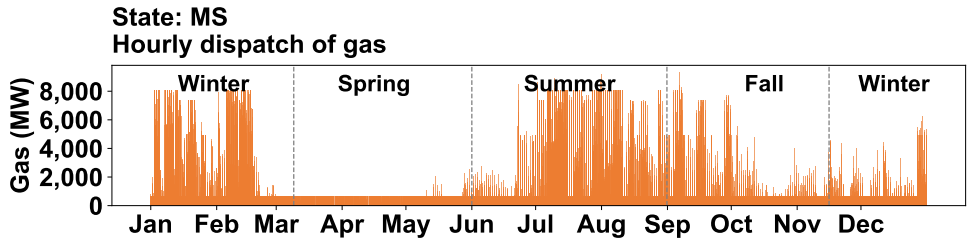

(c) Hourly dispatch of gas.

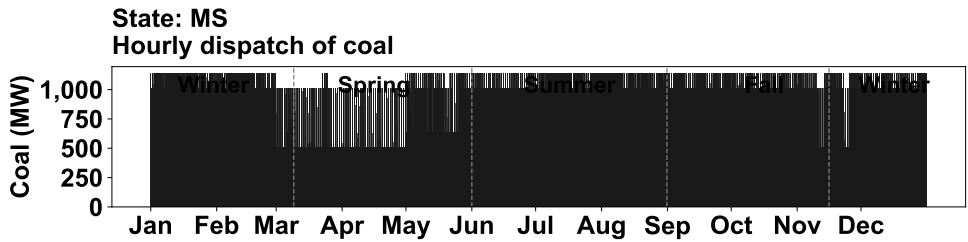

(d) Hourly dispatch of coal.

**Fig. 24**: Mississippi (MS).

28

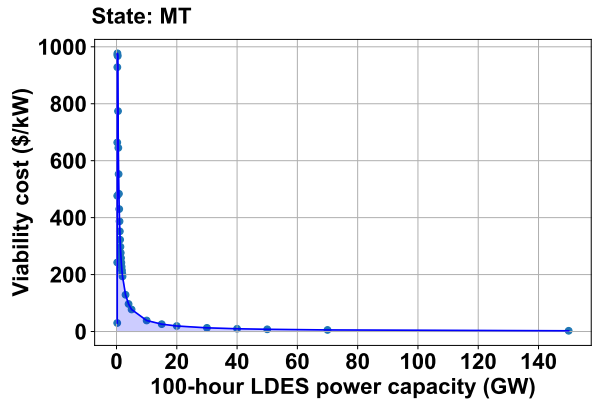

(a) Viability cost curve.

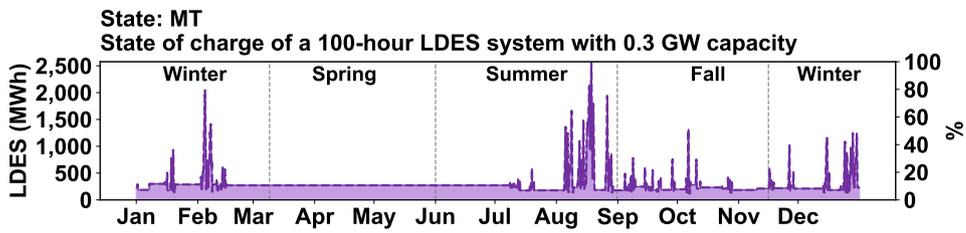

(b) State of charge.

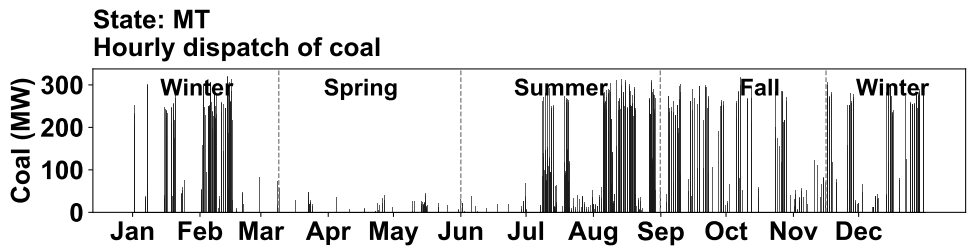

(c) Hourly dispatch of coal.

**Fig. 25**: Montana (MT).

29

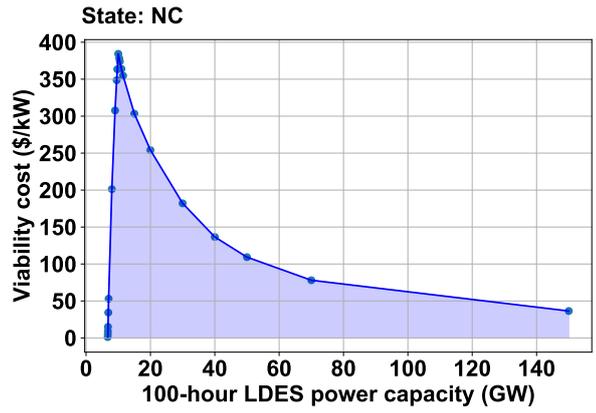

(a) Viability cost curve.

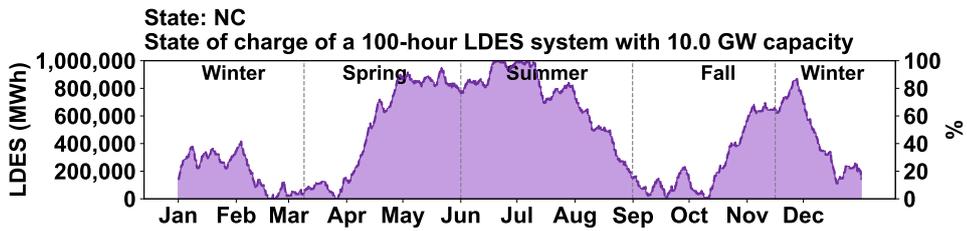

(b) State of charge.

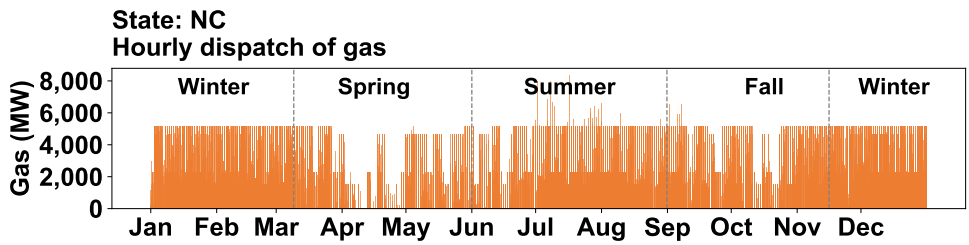

(c) Hourly dispatch of gas.

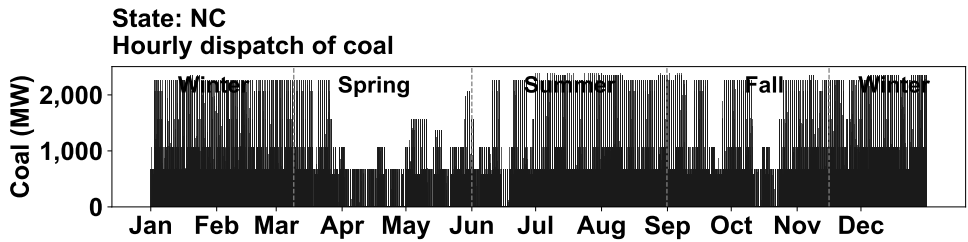

(d) Hourly dispatch of coal.

**Fig. 26**: North Carolina (NC).

30

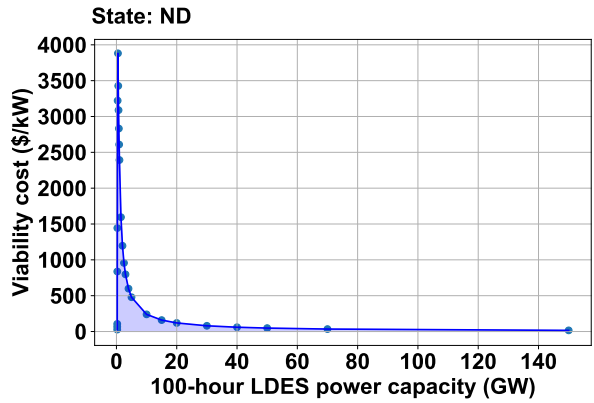

(a) Viability cost curve.

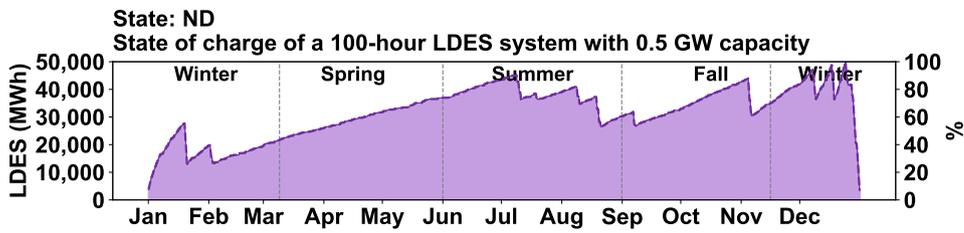

(b) State of charge.

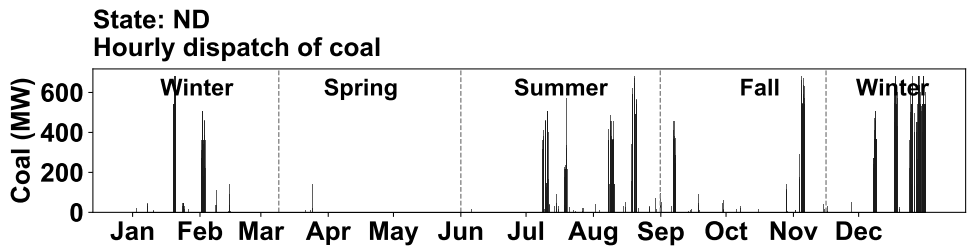

(c) Hourly dispatch of coal.

**Fig. 27**: North Dakota (ND).

31

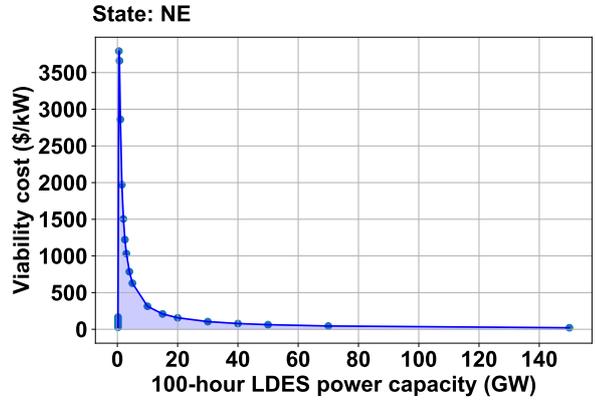

(a) Viability cost curve.

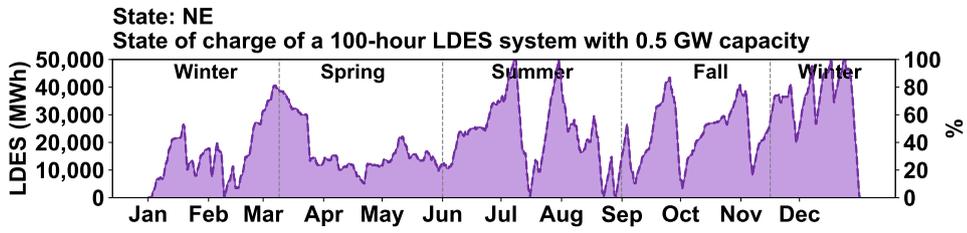

(b) State of charge.

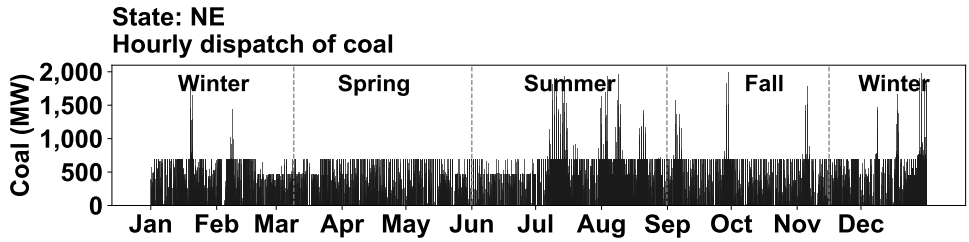

(c) Hourly dispatch of coal.

**Fig. 28**: Nebraska (NE).

32

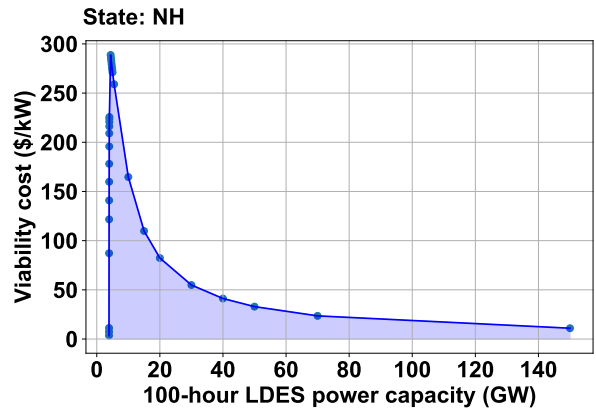

(a) Viability cost curve.

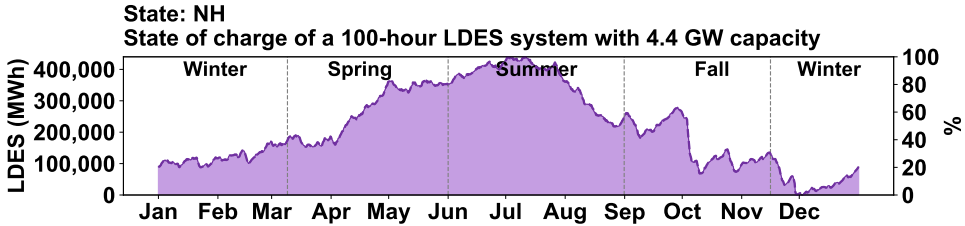

(b) State of charge.

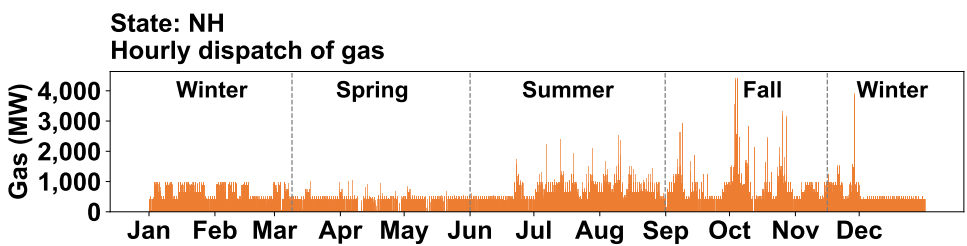

(c) Hourly dispatch of gas.

**Fig. 29**: New Hampshire (NH).

33

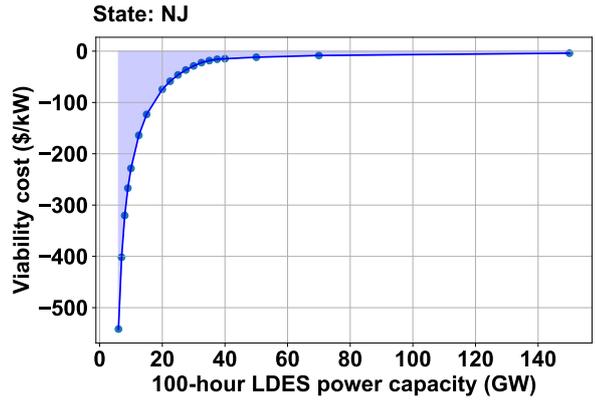

(a) Viability cost curve.

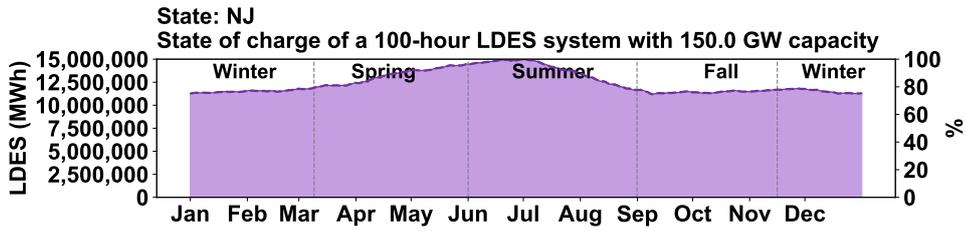

(b) State of charge.

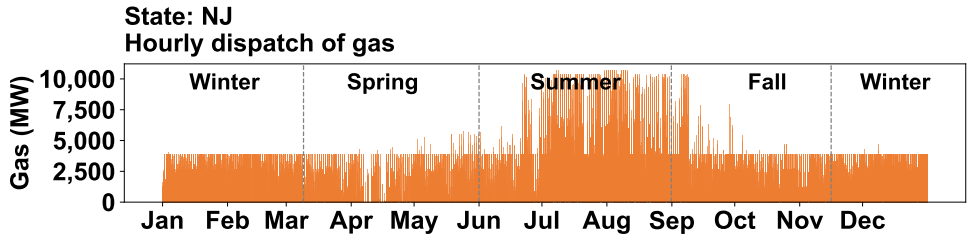

(c) Hourly dispatch of gas.

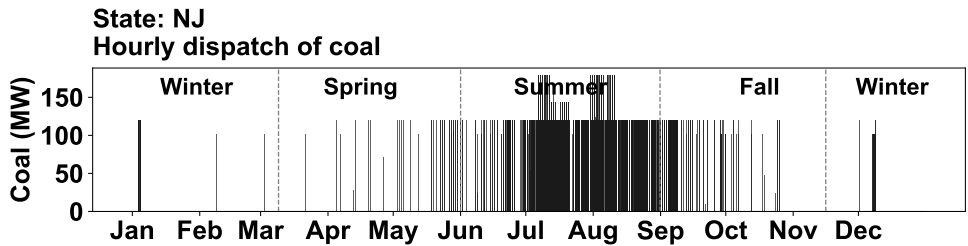

(d) Hourly dispatch of coal.

**Fig. 30**: New Jersey (NJ).

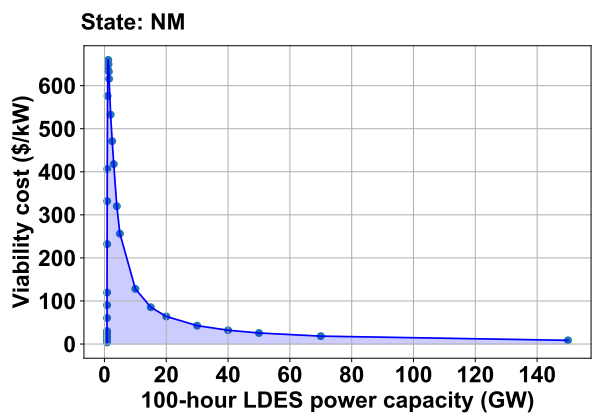
(a) Viability cost curve.

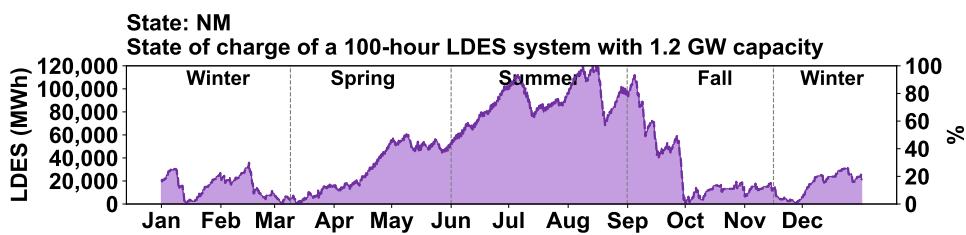
(b) State of charge.

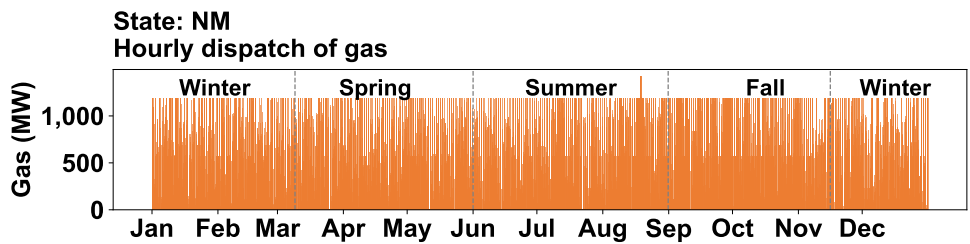
(c) Hourly dispatch of gas.

**Fig. 31**: New Mexico (NM).

35

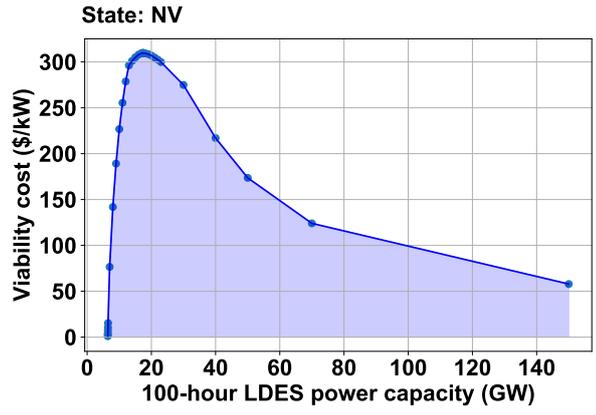

(a) Viability cost curve.

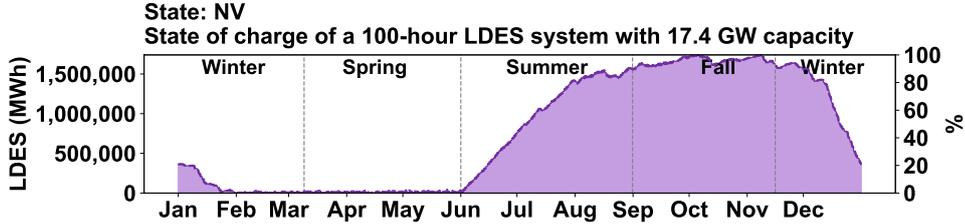

(b) State of charge.

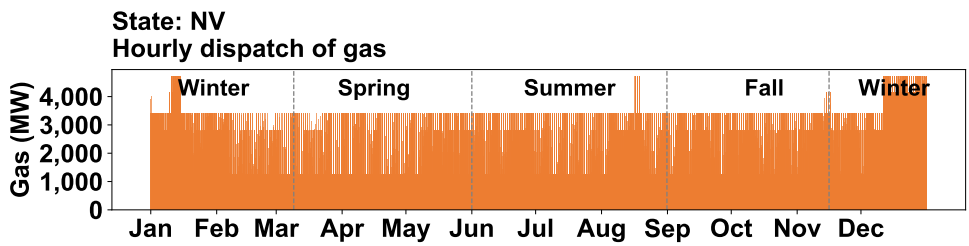

(c) Hourly dispatch of gas.

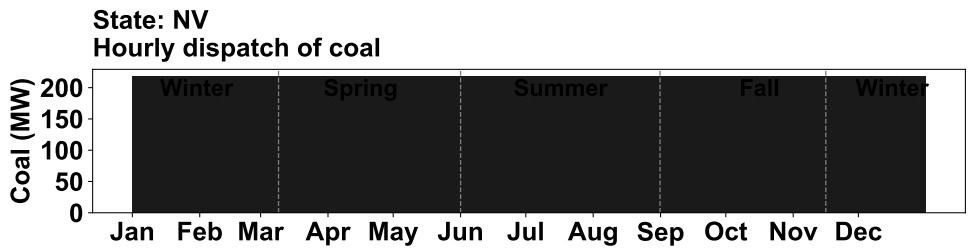

(d) Hourly dispatch of coal.

**Fig. 32**: Nevada (NV).

36

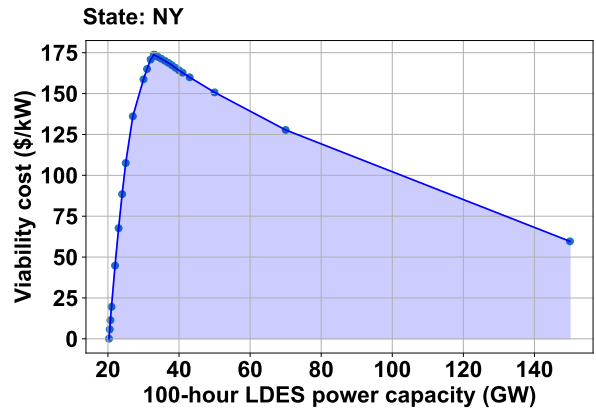

(a) Viability cost curve.

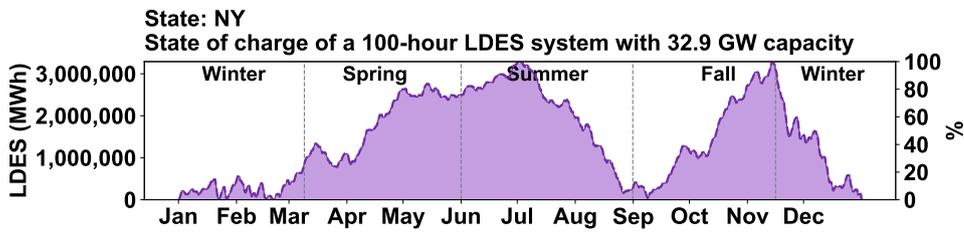

(b) State of charge.

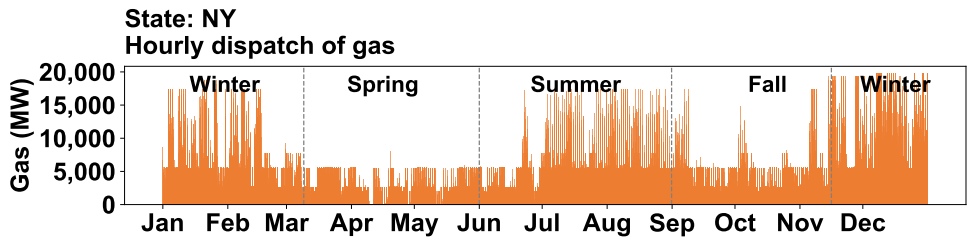

(c) Hourly dispatch of gas.

**Fig. 33**: New York (NY).

37

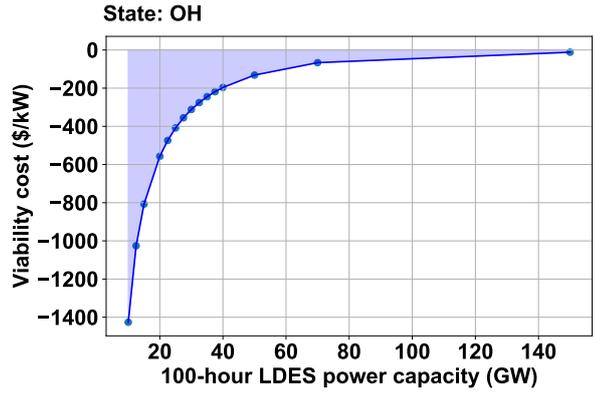

(a) Viability cost curve.

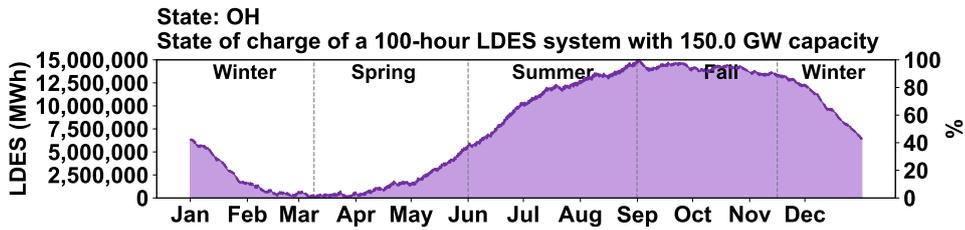

(b) State of charge.

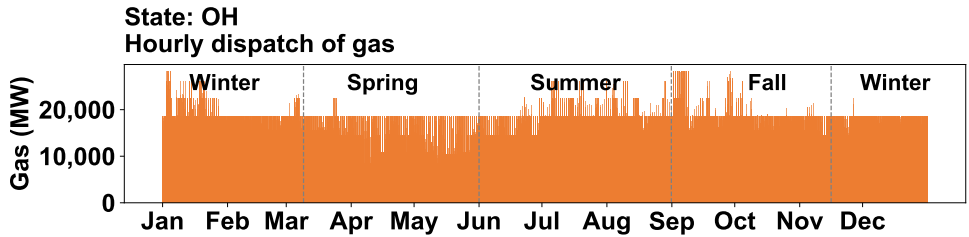

(c) Hourly dispatch of gas.

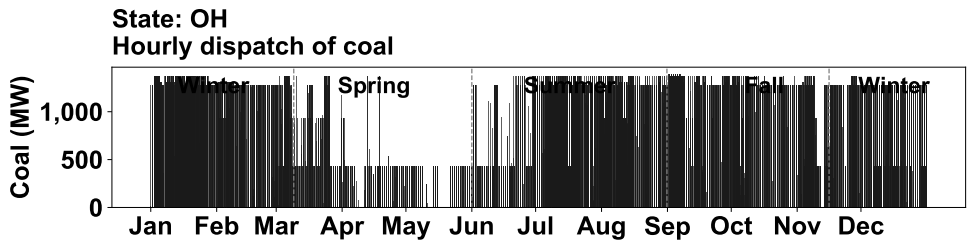

(d) Hourly dispatch of coal.

**Fig. 34**: Ohio (OH).

38

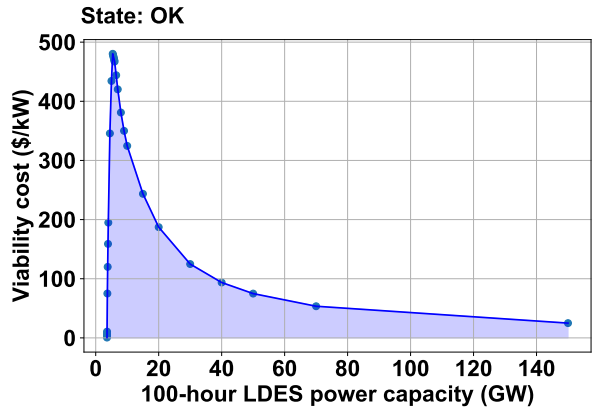

(a) Viability cost curve.

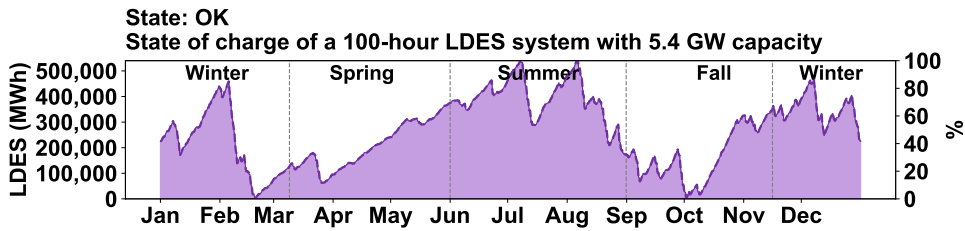

(b) State of charge.

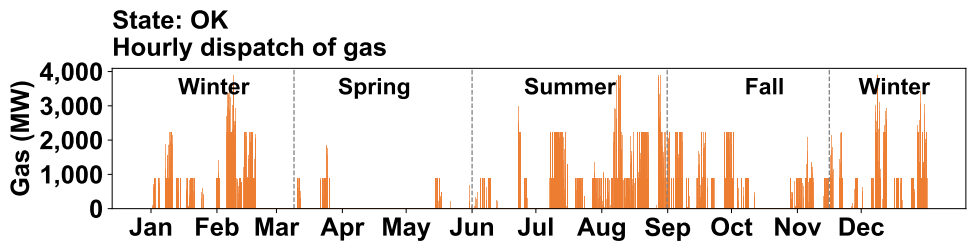

(c) Hourly dispatch of gas.

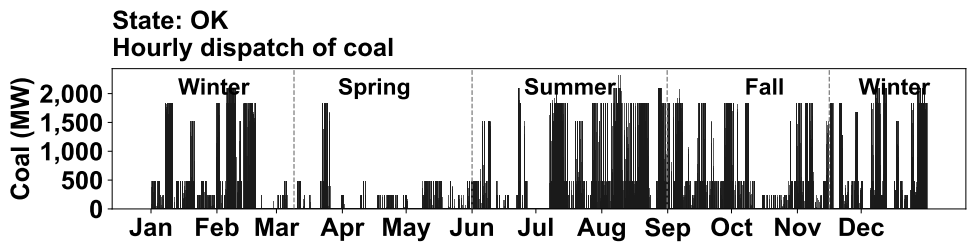

(d) Hourly dispatch of coal.

**Fig. 35**: Oklahoma (OK).

39

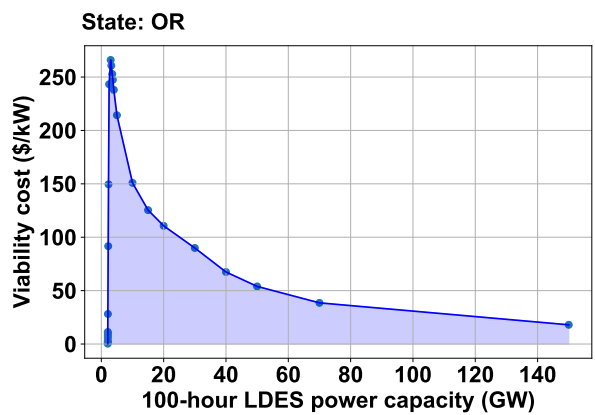

(a) Viability cost curve.

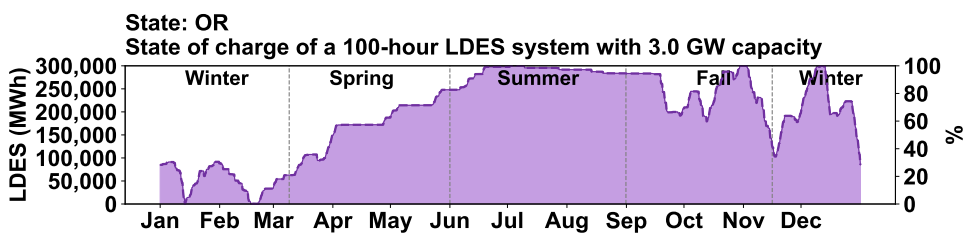

(b) State of charge.

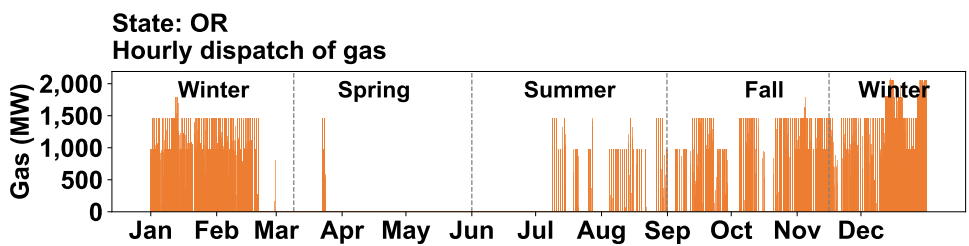

(c) Hourly dispatch of gas.

**Fig. 36**: Oregon (OR).

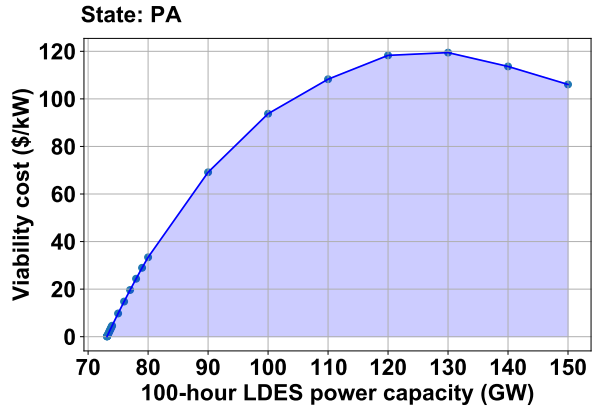

(a) Viability cost curve.

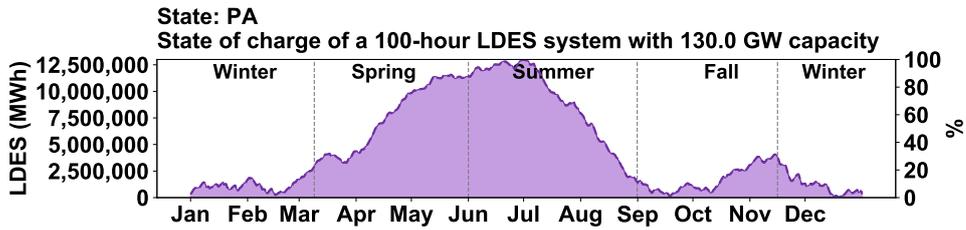

(b) State of charge.

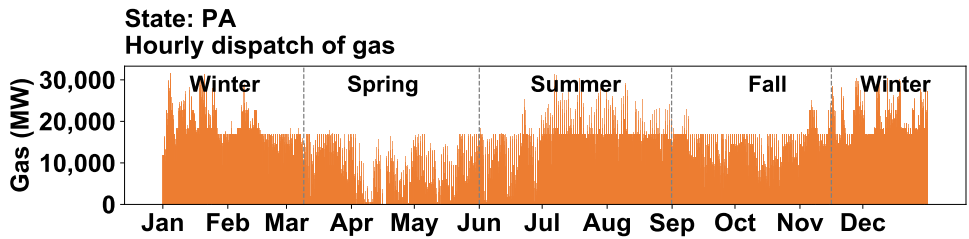

(c) Hourly dispatch of gas.

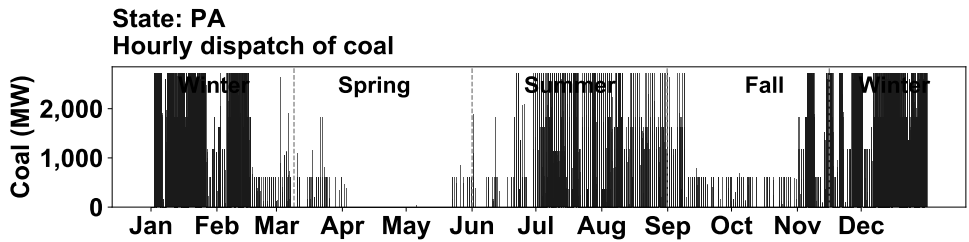

(d) Hourly dispatch of coal.

**Fig. 37**: Pennsylvania (PA).

41

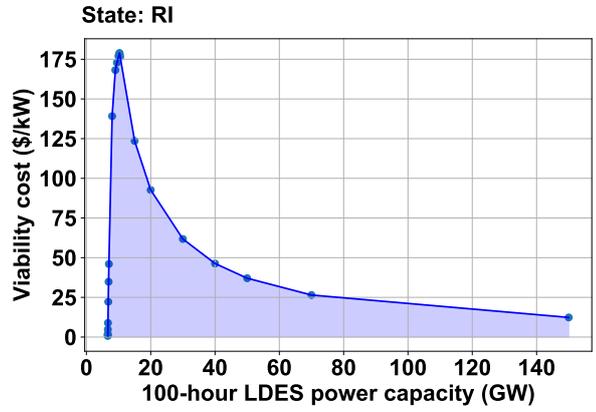

(a) Viability cost curve.

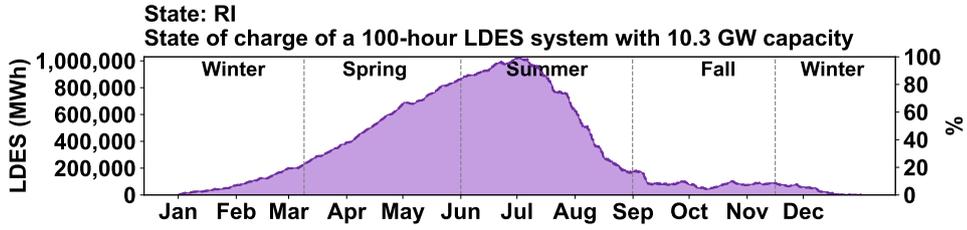

(b) State of charge.

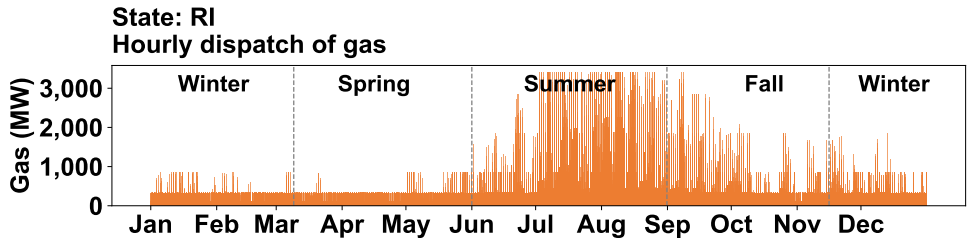

(c) Hourly dispatch of gas.

**Fig. 38**: Rhode Island (RI).

42

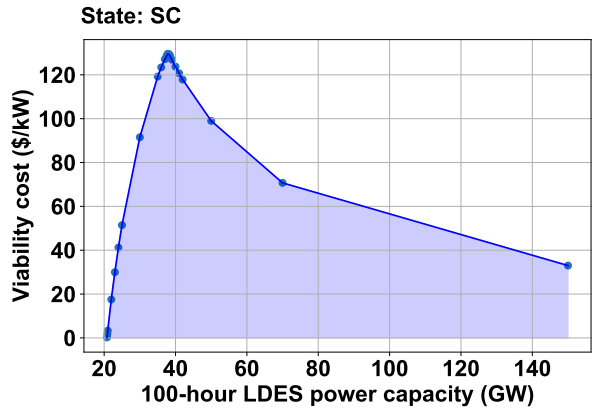

(a) Viability cost curve.

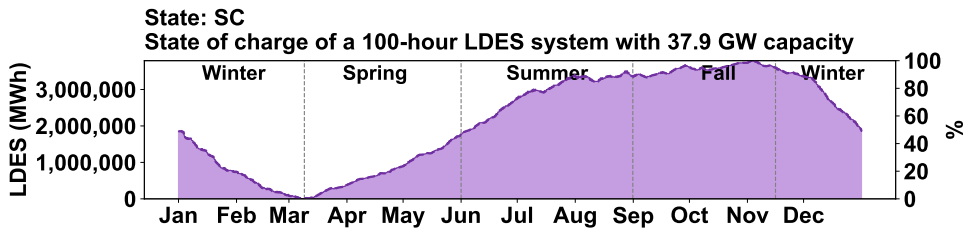

(b) State of charge.

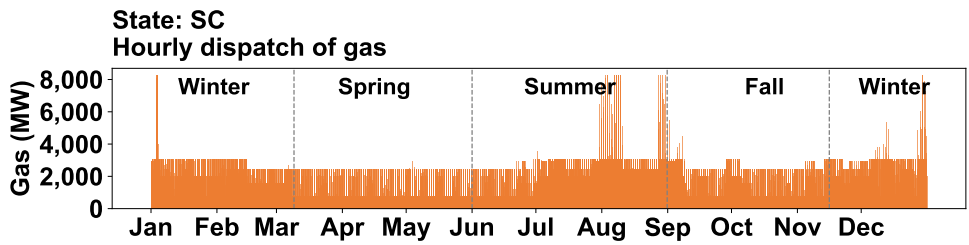

(c) Hourly dispatch of gas.

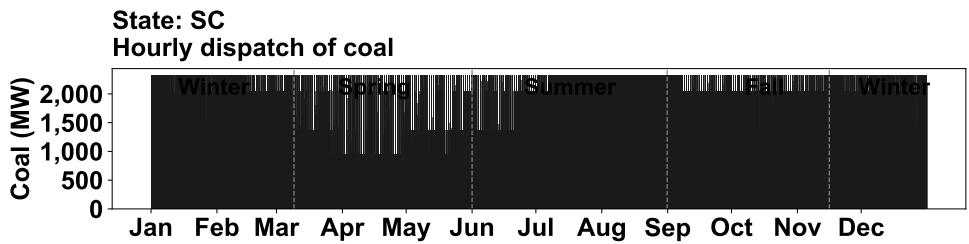

(d) Hourly dispatch of coal.

**Fig. 39**: South Carolina (SC).

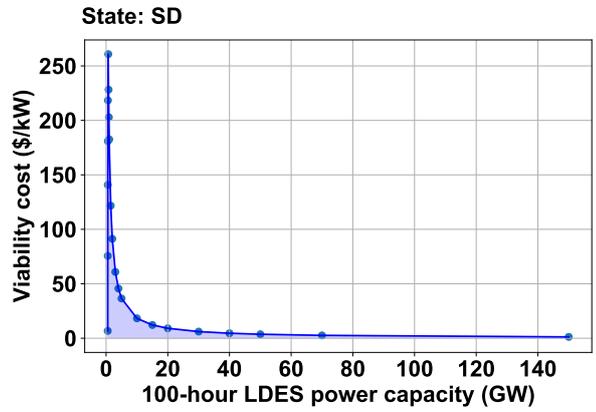

(a) Viability cost curve.

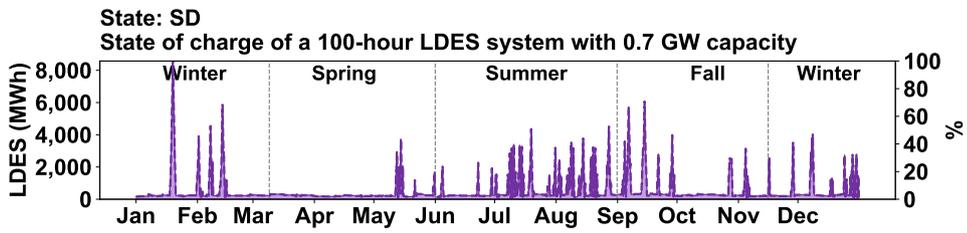

(b) State of charge.

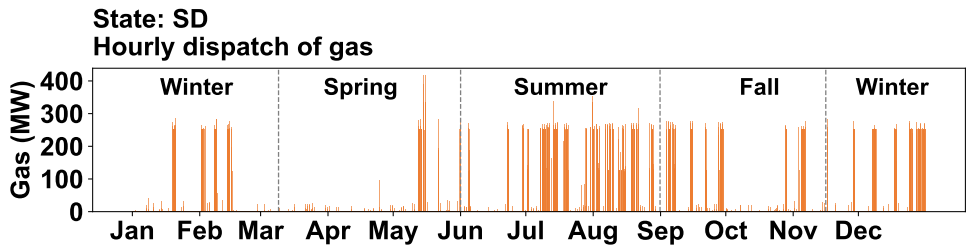

(c) Hourly dispatch of gas.

**Fig. 40**: South Dakota (SD).

44

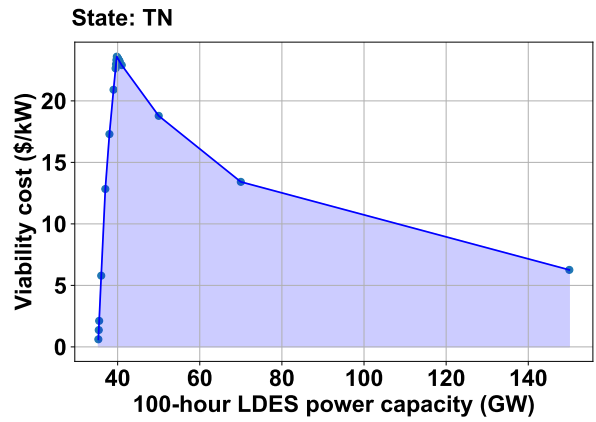

(a) Viability cost curve.

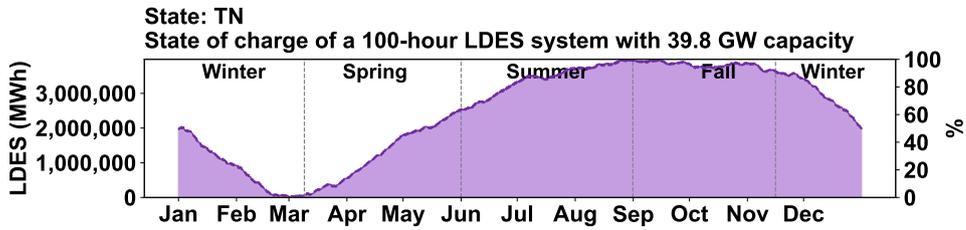

(b) State of charge.

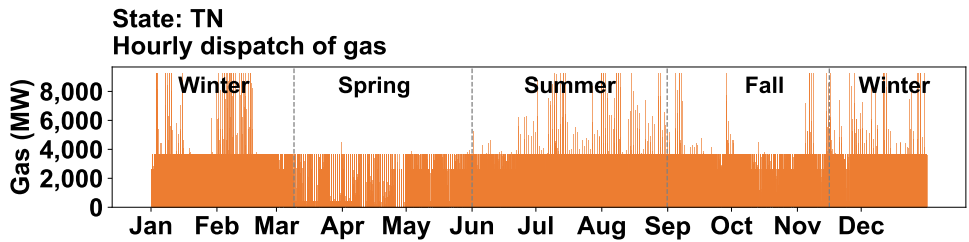

(c) Hourly dispatch of gas.

**Fig. 41**: Tennessee (TN).

45

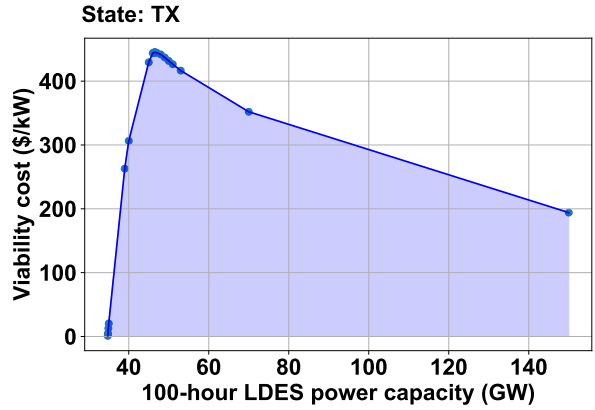

(a) Viability cost curve.

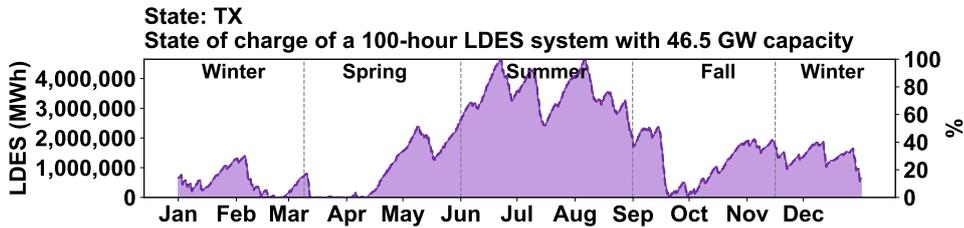

(b) State of charge.

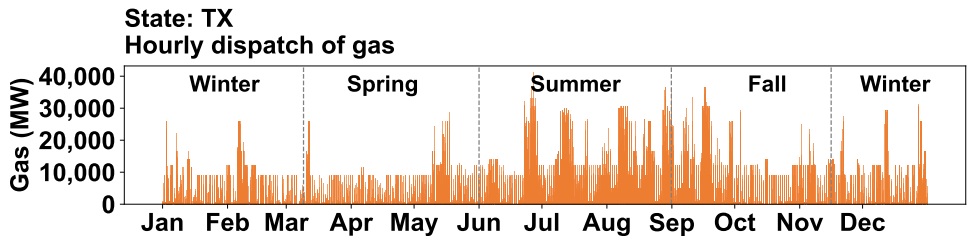

(c) Hourly dispatch of gas.

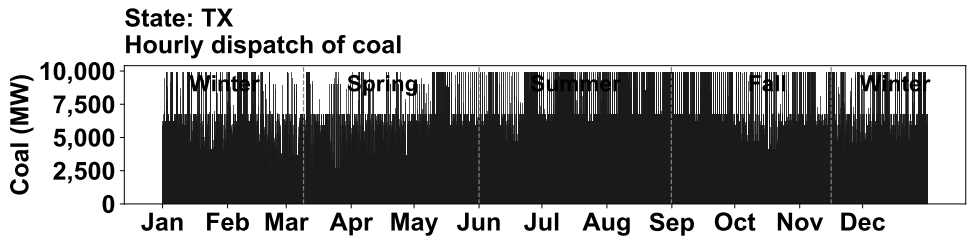

(d) Hourly dispatch of coal.

**Fig. 42**: Texas (TX).

46

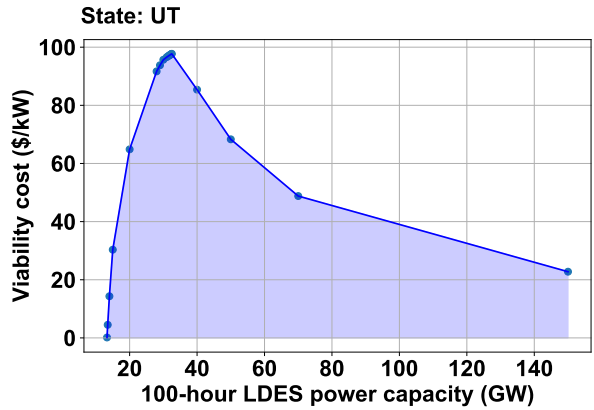

(a) Viability cost curve.

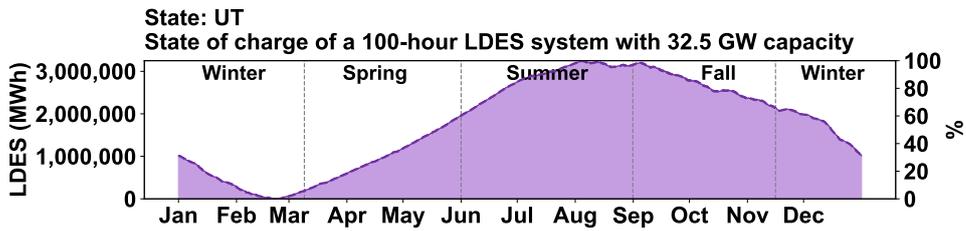

(b) State of charge.

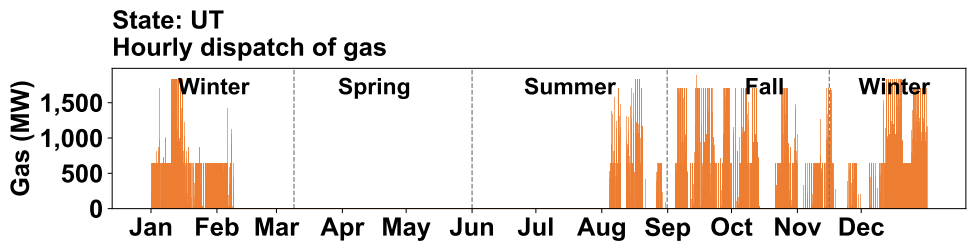

(c) Hourly dispatch of gas.

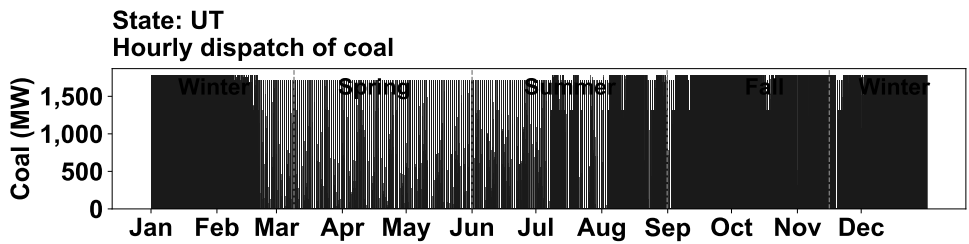

(d) Hourly dispatch of coal.

**Fig. 43**: Utah (UT).

47

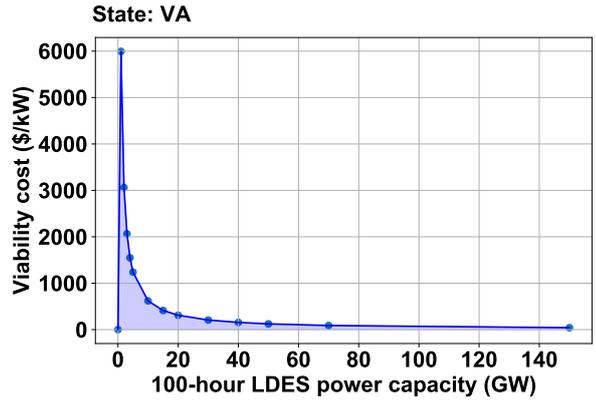

(a) Viability cost curve.

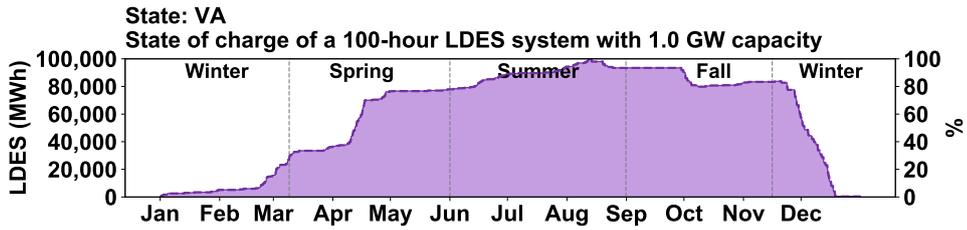

(b) State of charge.

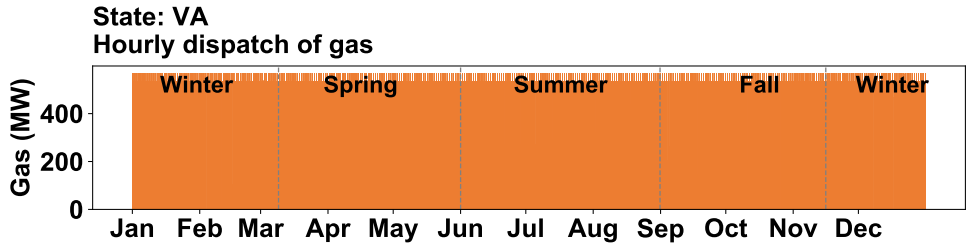

(c) Hourly dispatch of gas.

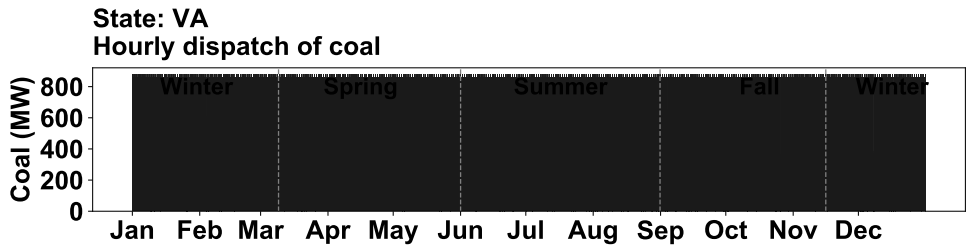

(d) Hourly dispatch of coal.

**Fig. 44**: Virginia (VA).

48

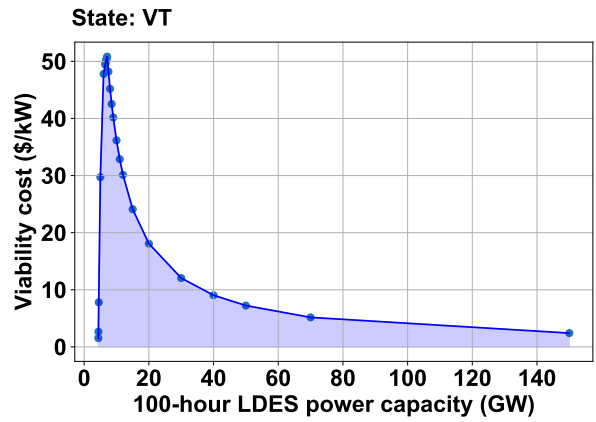

(a) Viability cost curve.

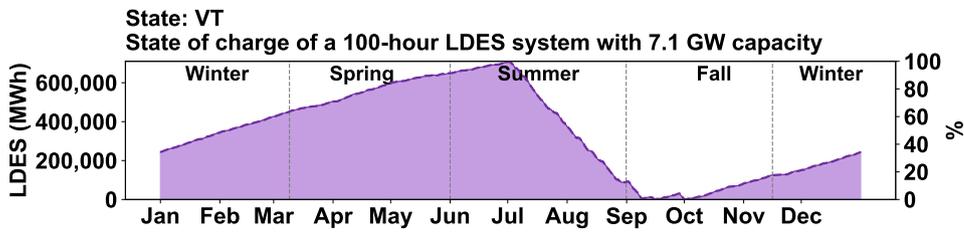

(b) State of charge.

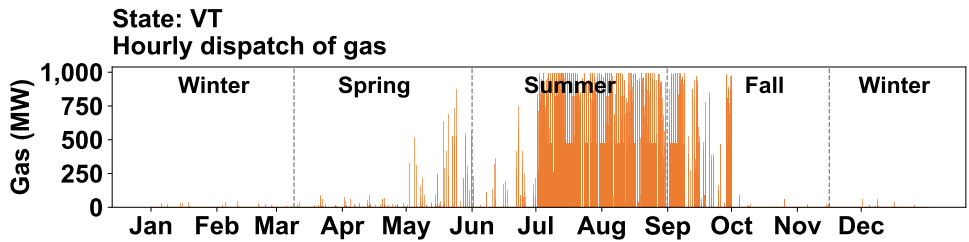

(c) Hourly dispatch of gas.

**Fig. 45**: Vermont (VT).

49

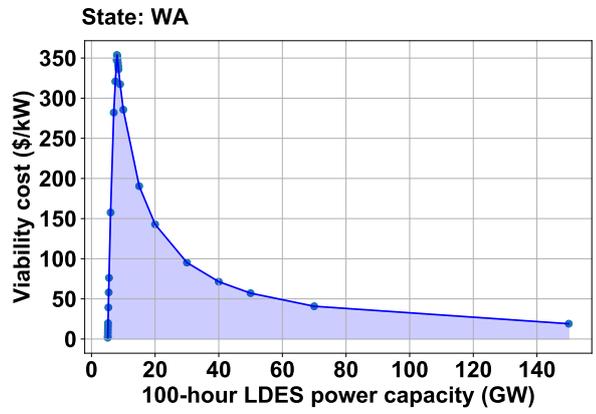

(a) Viability cost curve.

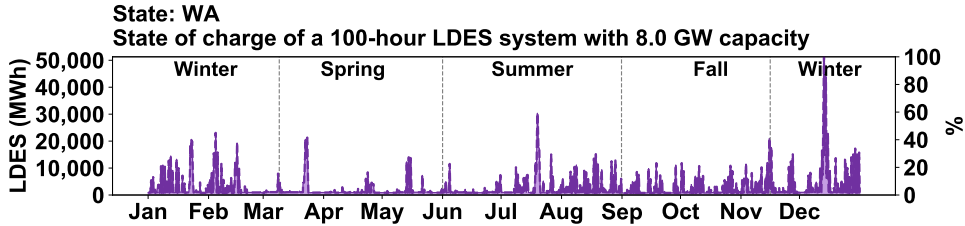

(b) State of charge.

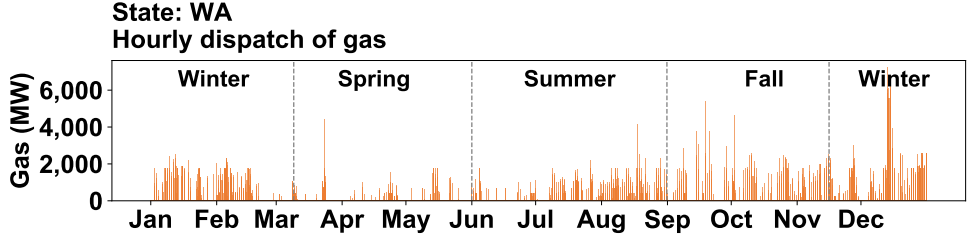

(c) Hourly dispatch of gas.

**Fig. 46**: Washington (WA).

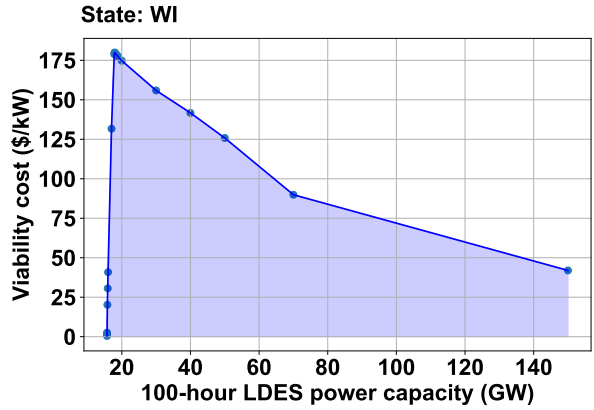

(a) Viability cost curve.

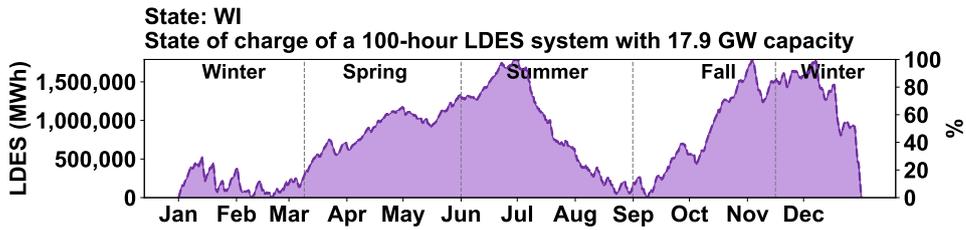

(b) State of charge.

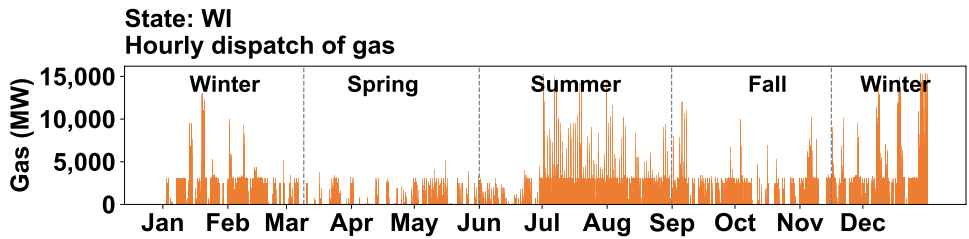

(c) Hourly dispatch of gas.

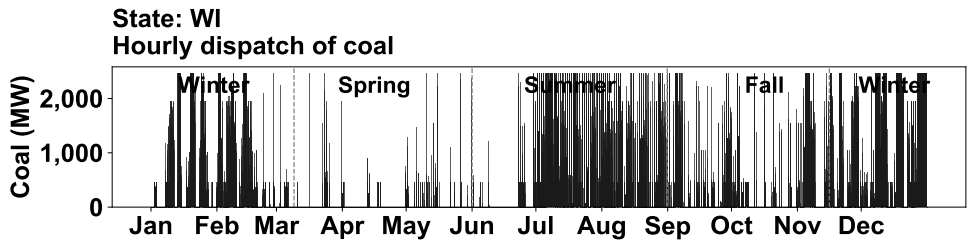

(d) Hourly dispatch of coal.

**Fig. 47**: Wisconsin (WI).

51

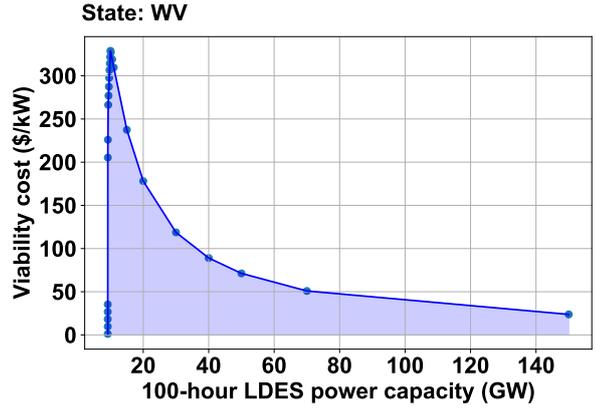

(a) Viability cost curve.

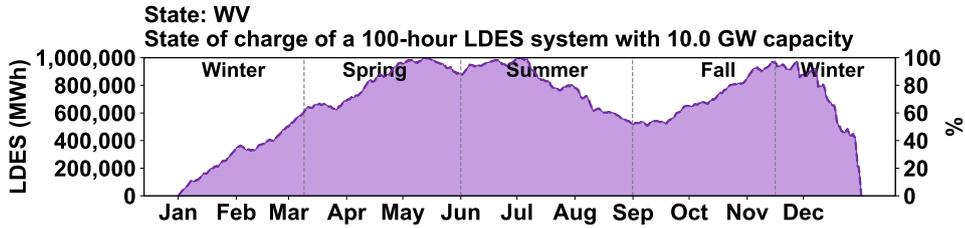

(b) State of charge.

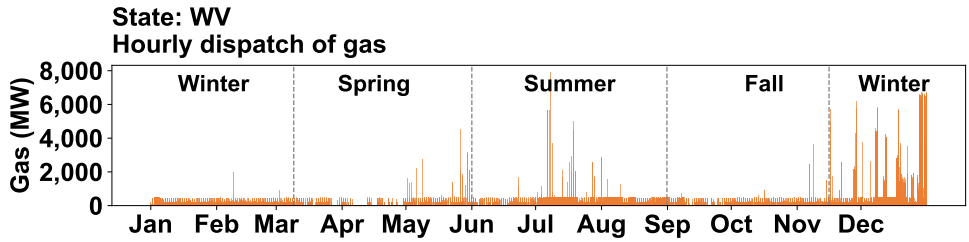

(c) Hourly dispatch of gas.

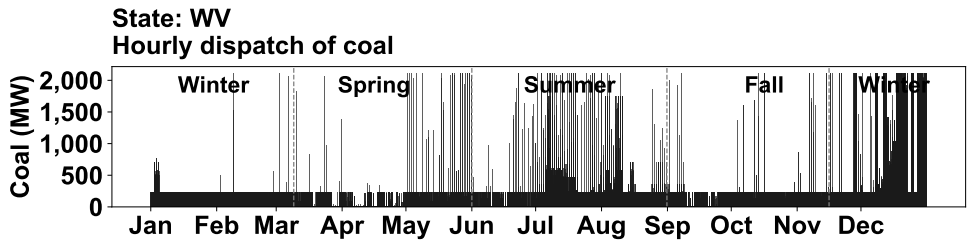

(d) Hourly dispatch of coal.

**Fig. 48**: West Virginia (WV).

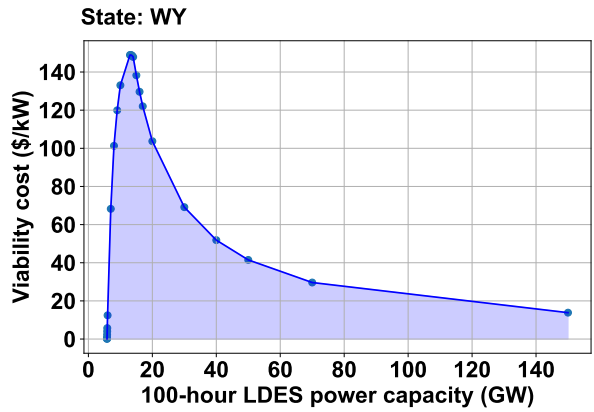

(a) Viability cost curve.

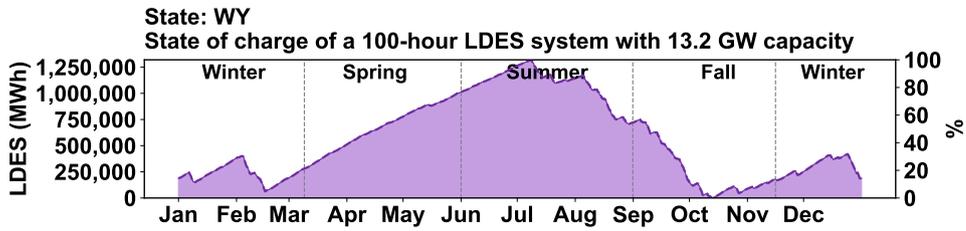

(b) State of charge.

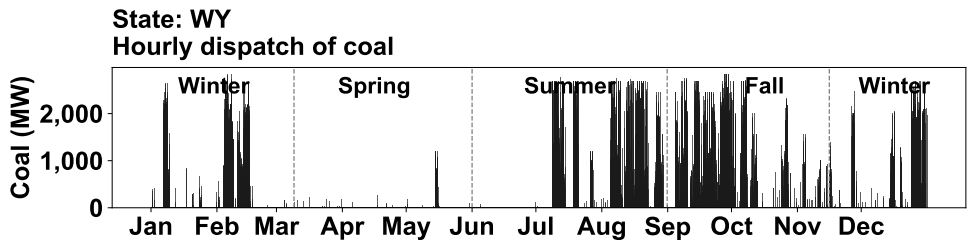

(c) Hourly dispatch of coal.

**Fig. 49**: Wyoming (WY).

53



\end{document}